\documentclass{amsart}    
\usepackage{amssymb}    
\usepackage{latexsym}    
    
\newcommand{\ZZ}{{\mathbb Z}}    
\newcommand{\RR}{{\mathbb R}}    
    
\newcommand{\NN}{{\mathbb N}}

\newtheorem{theorem}{Theorem}       
\newtheorem{remark}{Remark}       
\newtheorem{lemma}{Lemma}[section]       
\newtheorem{prop}[lemma]{Proposition}       
\newtheorem{coro}[lemma]{Corollary}       
\newtheorem{definition}[lemma]{Definition}       
\renewcommand\qedsymbol{$\Box$}    
    
\newcounter{smalllist}

\newcommand{\Xdot}{\dot{X}}
\newcommand{\Ydot}{\dot{Y}}

\newcommand{\Fgamma}{{ \mathcal{F} (\Gamma)}}
\newcommand{\Ogamma}{{ \mathcal{O} (\Gamma)}}
\newcommand{\Ggamma}{G(\Gamma) }
\newcommand{\Mgamma}{M(\Gamma) }

\newcommand{\Gugamma}{G_u (\Gamma) }
\newcommand{\Gmgamma}{G_m (\Gamma) }

\newcommand{\Gammanull}{\Gamma^{(0)}}
\newcommand{\Tau}{\mathcal{T}}

\newcommand{\graph}{{\bf  g}}

\newcommand{\Cl}{\mbox{Cl}}
\newcommand{\gros}{\mbox{$<$}}
\newcommand{\Igamma}{\mathcal{I }(\Gamma)}
\newcommand{\Ugamma}{\mathcal{V }(\Gamma)}
\newcommand{\Cred}{C_{\mbox{red}}^\ast}
\newcommand{\Gammatilde}{\widetilde{\Gamma}}

\newcommand{\diam}{\mbox{diam}}

\begin{document}    
\title[Construction of groupoids from inverse semigroups]{On an  order based construction of a groupoid from  an inverse semigroup }    
\author[ D.~Lenz]{Daniel H. Lenz}    
\maketitle    
\vspace{0.3cm}        
\noindent        
Fakult\"at f\"ur Mathematik, TU Chemnitz, D - 09107 Chemnitz, Germany
\\[2mm]  
E-mail: \mbox{dlenz@mathematik.tu-chemnitz.de }\\[3mm]        

\begin{abstract}    
We present a construction, which assigns two  groupoids,
$\Gugamma$ and $\Gmgamma$, to an inverse semigroup $\Gamma$. By
definition, $\Gmgamma$ is a subgroupoid (even a reduction)  of
$\Gugamma$.  

The construction unifies known constructions for
groupoids. More precisely,  the groupoid $\Gugamma$ is shown to be
isomorphic to the universal groupoid of $\Gamma$
introduced by Paterson. For $\Gamma$ arising from graphs resp. tilings,
the groupoid $\Gmgamma$ is the graph groupoid introduced by Kumjian
et al. resp. the   tiling groupoid introduced by Kellendonk. 

We obtain a  characterisation of open invariant sets in
$\Gmgamma^{(0)}$ in terms of certain order ideals of $\Gammanull$ for a large 
class of $\Gamma$ (including those arising from graphs and from tilings). If $\Gmgamma$ is essentially principal this gives a
characterization of the ideal structure of $\Cred(\Gmgamma)$ by a theory 
of Renault. In
particular, we then obtain necessary and sufficient conditions on
$\Gamma$ for simplicity of $\Cred(\Gmgamma)$. 

Our approach relies on a detailed analysis of the order structure of
$\Gamma$.

\end{abstract}    
    
\section{Introduction}   \label{Introduction}
This article is concerned with the construction of groupoids from
inverse semigroups and applications to graphs and tilings. 

The motivation for our study comes from  two sources. The first source is  work of Paterson \cite{Pat1,Pat2} and of Kellendonk
\cite{Kel,Kel2} (cf. \cite{KL, KL2} as well). Both Paterson and Kellendonk
present constructions assigning groupoids to inverse semigroups. The
relationship between their constructions is not quite apparent and it is
our aim  to present a unified view.
More precisely, we will show that the set $\Ogamma$ of directed subsets of
$\Gamma$  modulo an obvious equivalence relation is an inverse semigroup 
in the natural way (cf. \cite{Kel2} for related results as well).  Restricting the multiplication on $\Ogamma$ gives
then a groupoid $\Gugamma$ which is shown to be isomorphic to the universal
groupoid of Paterson.  It turns out that there is a subset of $\Ogamma$,
where the restriction of multiplication does not alter the multiplicable
pairs. This is the set of minimal elements of $\Ogamma$. Thus, this set
is a subgroupoid and even a reduction of $\Gugamma$. It will be   called
$\Gmgamma$. Under a rather mild assumption on  $\Gamma$ it can be shown
to agree with the groupoid $H_m (\Gamma)$ introduced by Kellendonk in
\cite{Kel2}.  

The second starting point for this work is given by recent
investigations on graphs and their associated groupoids by Kumjian et
al. \cite{KPRR}.  Kumjian et al. study (among other topics) the ideal
structure of $C^*$-algebras associated to graphs (cf. \cite{BPRS} and
\cite{HR} as well for different approaches to these algebras). An
important step in their analysis is a characterization of the open
invariant sets of the associated groupoid.

Here, it is our aim to present an inverse semigroup based approach to
these topics.

More generally, we present an abstract characterization of open
invariant sets of $\Gmgamma^{(0)}$ in terms of certain order ideals of
$\Gamma^{(0)}$ for quite general $\Gamma$. 

The connection to graphs is then made by showing that $\Gmgamma$ is
actually isomorphic to the graph groupoid of Kumjian et al. if $\Gamma$
is a suitable graph inverse semigroup. Restricting the  theory developed 
for general $\Gamma$ to the graph case,  we then recover the mentioned
results of Kumjian et al.. 

\medskip

Let us now discuss these points in some more detail. 
 As discussed by Renault \cite{Ren}, every topological groupoid $G$
 gives rise to an inverse semigroup viz its ample semigroup $G^a$. If the groupoid is ample, this ample semigroup determines the 
topology of the groupoid. This naturally raises the question whether
every semigroup can be faithfully represented as subsemigroup of  $G^a$
determining the topology of $G$ for a suitable groupoid $G$. 

This question has been thoroughly addressed by Paterson \cite{Pat1,Pat2}. It
turns out that such a representation is not unique in general and may
not even exist if $G$ is required to be Hausdorff \cite{Pat1}. On the
other hand, Paterson presents a  construction assigning what he 
calls the universal groupoid $G_u(\Gamma)$ to an inverse semigroup
$\Gamma$. This is an ample but not necessarily Hausdorff groupoid on
whose ample semigroup $\Gamma$ can be faithfully represented. It
determines all other so called $\Gamma$-groupoids and also determines
the representation theory of $\Gamma$ \cite{Pat2}. Paterson also shows
that $\Gugamma$ is Hausdorff if $\Gamma$ is $E$-unitary. 

This approach relies on a generalization of  Kumjians theory of
localization \cite{Kum} developed by Paterson in \cite{Pat1,Pat2}. Here,
a localization means a  suitable  action of an inverse semigroup
$\Gamma$ on a space $X$. As shown by Paterson this gives rise to an
$r$-discrete groupoid $G(X,\Gamma)$ whose $C^\ast$-algebra is isomorphic
to $C_0 (X)\times \Gamma$, where the partial crossed product is taken
the sense of Sieben \cite{Sie2}. In fact, results of Paterson
\cite{Pat2} and Sieben \cite{Sie2,Sie3} (cf. \cite{QS} as well) show,
roughly speaking,   that
there is a one-to-one  correspondence  between $r$-discrete groupoids
and (partial) actions of inverse semigroups on suitable spaces (at
least on the level of $C^\ast$-algebras).  Related topics  have also been
studied by Exel in \cite{Exe} and Nica in \cite{Nic}.

While $\Gugamma$ has very nice universal features its relation to
groupoids arising in concrete examples (as tiling groupoids or graph
groupoids) is not quite clear. 

On the other hand there is a different construction due to  Kellendonk
\cite{Kel2} assigning a groupoid $H_m (\Gamma)$ to an inverse
semigroup. In the context of tilings this construction gives the  tiling
groupoid.  Moreover, as we show below this construction also gives the
graph groupoid when applied to a suitable inverse semigroup associated
to a graph. 
The construction is based on suitable maximal ordered
sequences in $\Gamma$. 

Therefore, it is our first aim here to investigate the relation
between the groupoid introduced by Kellendonk and that introduced by
Paterson.  To do so we will have another and somewhat more systematic
look at the order based considerations done by Kellendonk in
\cite{Kel2}. This will give a unified construction to both groupoids and clarify their relationship. This
approach also  leads 
to interesting additional information on their respective topologies.

Our construction is  simple and does not use the theory of
localizations or any actions of $\Gamma$ on a space. It can rather be
seen as a kind of order completion of $\Gamma$ (cf. \cite{Bir} for
related material in general lattices). 
Moreover, it gives immediately  a lattice type condition on  $\Gamma$,  which we term  (L),  for $\Gugamma$ to be Hausdorff (Corollary \ref{Hausdorff}). Condition (L) just means suitable
existence of minima and is strictly  weaker than $E$-unitarity and $0$-$E$-unitarity. If
$\Gamma$ satisfies (L), the topology of $\Gugamma$ is not only Hausdorff
but even admits a different description. Namely, it can then simply
be described as the topology inherited from the product topology of
$\{0,1\}^\Gamma$ under the natural injection (Lemma \ref{einfach}).

If $\Gamma$ satisfies (L), then   the topology of
$\Gmgamma$ also  has  very nice  features and a particularly simple
basis  (Lemma \ref{wichtig}). This can be used to show that in this case 
the  groupoid $H_m (\Gamma)$ introduced by Kellendonk is isomorphic to
$\Gmgamma$.  In fact, $H_m (\Gamma)$ and $\Gmgamma$ always agree as sets 
but the topology might be different. These results clarify the  relationship of $H_m (\Gamma) $ and
$\Gugamma$. 

As (L) is satisfied for tiling inverse semigroups
our construction gives the tiling groupoid in this case. The condition
(L) is also satisfied for a suitable inverse semigroup associated to
graphs. In this case $\Gmgamma$ can be shown to be the graph groupoid
introduced by Kumjian et al.  (Theorem \ref{graphiso}). Thus, the
construction easily produces two important classes of groupoids. 

In the general case $\Gmgamma$ is a reduction of $\Gugamma$ (in the set
theoretical sense) on some invariant set $E$. We present two types of
general conditions  on $\Gamma$ for $E$ 
to be closed yielding that $\Gmgamma$ is a reduction of $\Gugamma$ in
the topological  sense (Proposition \ref{rad} and Proposition \ref{trap}). 
 Both types of  condition are   met in the tiling
case and in the graph case.  These general conditions are 
important because  they imply, in particular, that $\Gmgamma$ is locally
compact which   is not clear in general. 

Our study of $\Gmgamma$ leads us to a certain inverse subsemigroup
$\Gammatilde$ which is a quotient of $\Gamma$. This quotient
$\Gammatilde$ gives a precise version of how $\Gamma$ can be
considered as an inverse subsemigroup of the inverse semigroup of
$\Gmgamma$-sets in $\Gmgamma$.

We then investigate the lattice of open invariant sets in
$\Gmgamma^{(0)}$. This lattice plays a key role in the ideal theory of
the $C^\ast$-algebra $\Cred (\Gmgamma)$ due to a theory of Renault
\cite{Ren, Ren2}. Our main result there gives a characterization of
this lattice in terms of a lattice of certain order ideals in
$\Gamma^{(0)}$ (Lemma \ref{character}). This result can be used to
provide necessary and sufficient conditions on $\Gamma$ for
non-existence of non-trivial open invariant subsets of
$\Gmgamma^{(0)}$ (Lemma \ref{minimal}). If $\Gmgamma$ is essentially
principal, these Lemmata completely describe the ideal theory of
$\Cred (\Gmgamma)$ and give a necessary and sufficient condition on
$\Gamma$ for simplicity of $\Cred (\Gmgamma)$ (Theorem \ref{ideal}).

Unfortunately, we have not been able
to find convincing conditions on $\Gamma$ for $\Gmgamma$ to be
essentially principal.  For $\Gamma$  arising from (suitable)  graphs,
$\Gmgamma$ is essentially principal by  arguments of
\cite{KPRR} and we recover the results of \cite{KPRR}.  For
$\Gamma$ arising from tilings there is a simple condition for
principality of the corresponding groupoids. Thus, we can find a
description of ideals of $\Cred (\Gmgamma)$ in this case as well. 

\medskip

This paper is organized as follows: In Section \ref{Preliminaries}, we
review several known facts on (the order of)  inverse semigroups. In
particular, we show  that an inverse semigroup gives rise to two
groupoids, one arising by restricting the multiplication the other
consisting of minimal elements.

 Section \ref{Thebasic} contains the basic
constructions showing that the set of directed sets in an inverse
semigroup is again an inverse semigroup.  In Section \ref{ThegroupoidGu} 
we study certain aspects of $\Gugamma$ in some detail. This concerns in
particular the topology. We prove Lemma \ref{einfach} and  Corollary
\ref{Hausdorff}. Section \ref{neu} is devoted to a study of the inverse
semigroup $\Gammatilde$.
Section \ref{ThegroupoidGm} presents a general  study 
of $\Gmgamma$. It contains Lemma \ref{wichtig}, and Theorem \ref{consistent}. 

 Section \ref{Applicationgraph} is devoted to applications to
 graphs. We prove Theorem \ref{graphiso} and show how the material of
 the preceeding sections can be used to recover some results of
 Kumjian et al..  Finally, in Section \ref{tilings}, we recall results
 of \cite{Kel,Kel2} on tilings and provide a study of ideal theory of
 the algebras $\Cred (\Gmgamma)$ in this case. This underlines the
 similarity between the tiling case and the graph case.

\section{Preliminaries} \label{Preliminaries}
In this section we fix some notation and review  basic facts from 
inverse semigroup theory. In particular, we collect  properties
concerning the order structure of an inverse semigroup. For the general
theory of inverse semigroups  we refer the reader to  e.g. \cite{ Law, Pat2}. There, one
can also find those proofs which are omitted below.  The results of this
section will be used tacitly in the sequel. 

\medskip

Let $\Gamma$ be an inverse semigroup. This means that $\Gamma$ is a
semigroup and to each $x\in \Gamma$ there exists a unique $i(x)\in \Gamma$ with $x\, 
i(x) \,x =x$ and $i(x)\, x \, i(x) = i(x)$. The element $i(x)$ is denoted by
$x^{-1}$ and called the inverse of $x$. The map $x\mapsto x^{-1}$ is an
involution. By $\Gammanull$ we denote the units of $\Gamma$, i.e. the
set of $p$ with $p= p p^{-1}$. The units are just the idempotents and
commute. On $\Gamma$ we have the relation $\prec$, where $x\prec y$,
whenever $x y^{-1} = x x^{-1}$. If $x\prec y$, then $x$ is said to be a
precessor of $y$ and $y$ is said to be a successor of
$x$. Alternatively, $x$ is said to be smaller than $y$. 
The following proposition is well known
and easy to prove.

\begin{prop}\label{ordeins} For $x,y\in \Gamma$ the following are equivalent:\\
(i) $x\prec y$. (ii) $x^{-1} \prec y^{-1}$. (iii) $\exists p\in
\Gammanull$  $x =p y$. (iv)  $\exists q\in \Gammanull$  $x =y q$.
\end{prop}
Using this proposition, it is not hard to see that $\prec $ is an order
on $\Gamma$, i.e. a reflexive, transitive relation s.t. $x\prec y $ and
$y\prec x$ implies $x=y$. Moreover, we easily deduce from the
proposition the following proposition.
\begin{prop}\label{ordzwei}If  $x_1 \prec y_1$ and $x_2\prec y_2$, then 
  $x_1 x_2 \prec y_1 y_2$. 
\end{prop}
As $\prec$ is an order, there can not be more than one $z\in \Gamma$
with $z\prec x$ for every $x\in \Gamma$. A simple calculation shows that 
such a $z$ satisfies $z x=z 
= x z$ for every $x\in \Gamma$. Therefore, it  will be denoted by
$0$. $\Gamma$ is said to be an inverse semigroup with zero if it
contains such a $z$. In the sequel we will sometimes write conditions of 
the form $0\neq z\in \Gamma$. This is meant to mean that $z$ is not zero 
if $\Gamma$ has a zero and to be a vacuous condition if $\Gamma$ does
not contain a zero. 

An element $x\in \Gamma$ is called minimal if it is not zero and  $y\prec x$ and $y\neq 0$
implies $y=x$. The set of minimal elements in $\Gamma $ is denoted by
$\Gamma_{{\rm min}}$. The following proposition follows easily from the
above two propositions.
\begin{prop}\label{ordeinseins} For $x\in \Gamma$ the following are equivalent:\\
(i) $x$ is minimal. (ii) $x^{-1}$ is minimal. (iii) $x^{-1} x$ is
minimal. (iv) $x x^{-1}$ is minimal.
\end{prop}

Moreover, we have the following result.

\begin{prop}\label{useful} For $x, y\in \Gamma_{\rm min}$ the following are
  equivalent:\\
(i) $ x y \neq 0$. (ii) $x^{-1} x = y  y^{-1}$. (iii) $x y\in
\Gamma_{\rm min}$. 
\end{prop}
{ \it Proof.} (i) $\Longrightarrow$ (ii). By $0\neq x y= x y y^{-1} y$,
we have $x y y^{-1} \neq 0$. By minimality of $x$, this gives $x = x
y^{-1} y$, which in turn implies $x^{-1} x = x^{-1} x y
y^{-1}$. Similarly, we calculate $y y^{-1} =x^{-1} x y
y^{-1}$ and (ii) follows.\\
(ii)$\Longrightarrow$ (iii). By (ii), we have $x= x x^{-1} x= x y
y^{-1}$ and  $x y\neq 0$ follows. Moreover, by minimality of $y$ we have 
$x y p= x y$ for every $p\in \Gammanull$ with $x y p\neq 0$. \\
(iii)$\Longrightarrow$ (i). This is immediate from the definition of
minimality. \hfill \qedsymbol. 

\medskip

The order $\prec$ will in general not be a semilattice. However, one can 
still ask for the existence of a largest common precessor of $x$ and $y$ 
given that there exist common precessor. If such a largest common precessor  exists it must be unique
and will be denoted by $x\wedge y$. The existence of such largest
precessors will be of crucial importance in our dealing with topological 
properties of the groupoids in question. We include a short
discussion. 

\begin{definition} An inverse semigroup is said to satisfy the lattice  condition (L), if for any $x,y\in \Gamma$, with a common precessor not equal to
  zero there exists a largest common precessor.
\end{definition}

The following definition gives a well known criterion for the existence
of largest common precessors. Recall that an order ideal $\mathcal{I}$
in $\Gamma$  is a set with $\{y: y\prec x \}\subset \mathcal{I}$ for
every $x\in \mathcal{I}$.  By an inverse order ideal, we mean a set $\mathcal{I}$
in $\Gamma$ with $\{y: x\prec y \}\subset \mathcal{I}$ for
every $x\in \mathcal{I}$.

\begin{definition} An inverse semigroup  called $E$-unitary if
  $\Gammanull$ is an inverse order ideal. An inverse semigroup is called
  $0$-$E$-unitary if $\Gammanull\setminus\{0\}$ is an inverse order ideal. 
 \end{definition}
\begin{remark}{\rm  (a) Apparently, we have $E$-unitary $\Longrightarrow $
    $0$-$E$-unitary. But, of course,  $E$-unitary is essentially used
    for inverse semigroups without zero, while $0$-$E$-unitarity is the
    right notion for inverse semigroups with zero.\\
(b) The condition  (L) is strictly weaker than
$E$-unitarity and $0$-$E$-unitarity. This can be seen by considering a groupoid $G$ and
inverse    subsemigroups of  the inverse  semigroup $S(G)$ of its $G$-sets }
\end{remark}
\begin{prop} Let $\Gamma$ be a E-unitary (or $0$-E-unitary).  If $x,y\in \Gamma$ have a
  common precessor $z$ not equal to zero, then 
    there exists a largest such $z$. It is given
  by $y x^{-1} y= x y^{-1} x$.
\end{prop}

While E-unitarity resp. $0$-$E$-unitarity has been used, when studying topological properties
of groupoids associated to inverse semigroups \cite{Kel2,Pat2}, it turns
out that   our considerations  need only the  weaker
condition (L).

Actually, the existence of  $x \wedge y$  can be shown for arbitrary $\Gamma$ under
certain conditions on $x,y$. This is investigated in the next proposition. 
 \begin{prop}\label{minimum} Let $x,y\in \Gamma$  with a common
  successor $z\in \Gamma$ be given. If there exists a common precessor
  of $x$ and $y$ not equal to zero then there exists a largest common
  precessor. It is given by $x x^{-1} yy^{-1} x= x x^{-1} y y^{-1} y$.
\end{prop}
{\bf Proof.} As $x$ and $y$ are smaller than $z$ the elements $p=x
x^{-1},q= y y^{-1}$ belong to  $\Gammanull$ with $x=pz $ and $y=q z$. As $x$ and
$y$ have common precessors not equal to zero, the product $p q z$ is not
zero. It is obviously smaller than $x$ and $y$. Moreover, it is not hard
to show that any $r$ smaller than $x$ and $y$ must be smaller than $p q
z$.\hfill \qedsymbol

\medskip


Let us now turn to groupoids. Recall that a groupoid is a set $G$
together with a partially defined associative    
multiplication $*$ and an
involution $x\mapsto x^{-1}$ satisfying the following conditions
\cite{Ren}: 

\begin{itemize}
\item[(G1)] $(x^{-1})^{-1}=x$.
\item[(G2)]  If $x*y$ and $y*z$ exist, then $x*y*z$ exists as well.
\item[(G3)]  $x^{-1}* x$ exists and if $x* y$ exists as well then $x^{-1}*x* 
  y=y$.
\item[(G4)] $x* x^{-1}$ exists and if $z*x$ exists as well then $z* x* x^{-1}=z$.
\end{itemize}

Now, there is an immediate way to construct two groupoids from
$\Gamma$. We start with the following result contained in e.g.     Proposition 1.0.1 of \cite{Pat2}. 

\begin{prop}  Let $\Gamma$ be an inverse semigroup. Then, $\Gamma$ with
  its usual inversion and multiplication defined by $x * y= x y$ if and
  only if $x^{-1} x= y y^{-1}$ is a groupoid.
\end{prop}

\begin{definition}\label{defeins} Let $\Ggamma$ be the groupoid associated to
  $\Gamma$ in the foregoing proposition. 
\end{definition}

\begin{prop} The set of minimal elements of $\Gamma$ with involution
  from $\Gamma$ and multiplication defined  whenever the product is not
  zero,  is a subgroupoid of 
  $\Ggamma$.
\end{prop}
{\it Proof.} By Proposition \ref{ordeinseins} the set of minimal elements is 
closed under inversion. By Proposition \ref{useful}, it is further
closed under the multiplication in $\Ggamma$. \hfill \qedsymbol

\begin{definition} Let $\Mgamma$ be the subgroupoid of $\Ggamma$
  consisting of minimal elements of $\Gamma$.  
\end{definition}

To state our next result, we recall some more facts. A subset $E$ of the 
units $G^{(0)}$ of  a groupoid $G$ is called invariant if for $e\in E $ and $g\in G$
with $e = g^{-1} g$ the element $g g^{-1}$ belong to $E$ as well. If $E$ 
is invariant, the (set theoretic) reduction $G_E$ of $G$ to $E$ is the subgroupoid of
$G$ consisting of all elements $g\in G$ with $g^{-1} g\in E$ (wich by
invariance implies $g
g^{-1} \in E$ as well). In the context of topological groupoids, the invariant
set $E$ is further required to be closed in $G^{(0)}$. We will then
speak of the topological reduction $G_E$. 

\begin{prop}\label{reduction} With $E\equiv \Gammanull \cap \Gamma_{\rm
    min}\equiv \Gammanull_{\min}$ the equality  $\Mgamma= \Ggamma_{E}$ holds.
\end{prop}
{\it Proof.} The set $E
$ is invariant by
Proposition \ref{ordeinseins}. Thus, $\Ggamma_{E}$ is indeed a groupoid. As the groupoid structure is induced
by the groupoid structure of $\Ggamma$ for both $\Mgamma$ and  $\Ggamma_{E}$, it is enough to show that the underlying sets are
equal. This is easy. \hfill \qedsymbol

\medskip

In a groupoid $G$, the set of its $G$-sets is an inverse semigroup. This 
inverse semigroup will be denoted by $S(G)$.

\section{The basic construction} \label{Thebasic}
In this section we will show that the set of downward directed subsets
of $\Gamma$ modulo a certain equivalence relation is an inverse
semigroup with respect to the obvious multiplicative structure.  Lemma
\ref{basic} is strongly related to results of \cite{Kel}. This is
discussed at the end of Section \ref{ThegroupoidGm}.

\medskip

\begin{definition} A subset $A$ in $\Gamma$ is called (downward)
  directed if for any $x,y\in A$, there exists a $z\in A$ with $z\prec
  x,y$. The set of all directed subsets of $\Gamma$ is denoted by $\Fgamma$.
\end{definition}

On $\Fgamma$, we define the relation $\prec$ by $A\prec B$, if for any
$b\in B$, there exists an $a\in A$ with $a\prec b$. Moreover, we define 
$A B$ by $AB\equiv \{a b: a\in A, b\in B\}$ and $A^{-1}$ by
$A^{-1}=\{a^{-1}: a\in A\}$. The  corresponding sets are indeed directed
  by the results the first section. Moreover, we set $A\sim B$, whenever 
  $A\prec B$ and $B\prec A$. It is not hard to see that $\sim$ is an
  equivalence relation on $\Fgamma$. We set $\Ogamma=\Fgamma/\sim$.
  Representatives of $X, Y\in \Ogamma$ will be denoted by $\Xdot$ and
  $\Ydot$. The class of $A\in \Fgamma$ will be denoted by $[A]$. 
On
  $\Ogamma$, we define a multiplication by 
$$ X Y\equiv [ \Xdot \Ydot],$$
where $\Xdot$ and $\Ydot$ are arbitrary representatives of $X$ and
$Y$. It is easy to check that this a a well defined associative
multiplication. Moreover, we define a map $i: \Ogamma\longrightarrow
\Ogamma$ by $i(X)\equiv [\Xdot^{-1}]$, where again $\Xdot$ is a
representative of $X$ and this is well defined. The following lemma
holds.

\begin{lemma} \label{basic} The set $\Ogamma$ with multiplication and inversion
  $X^{-1} \equiv i(X)$ is an inverse semigroup. The relation $X\prec Y$
  holds for $X,Y\in \Ogamma$ if $\Xdot \prec \Ydot $ holds for some (all)
  representatives $\Xdot$ of $X$ and $\Ydot$ of $Y$. 
\end{lemma}
{\it Proof.}
 We first show that each $X\in\Ogamma$ has a unique inverse
given by $i(X)$. Existence follows easily from 
$$\Xdot =\{x x ^{-1} x : x\in \Xdot \}\sim \{x_1 x_2^{-1} x_3: x_1, x_2, 
x_3\in \Xdot\} = \Xdot \Xdot^{-1} \Xdot.$$
To show uniqueness, let $X$ and $Y$ be given with representatives $\Xdot 
$ and $\Ydot$ and assume  (1)  $X Y X= X$  and (2) $Y X Y =Y$. By (1), we
have
$$ \Xdot^{-1} \sim \Xdot^{-1} \Xdot \Xdot^{-1} \sim  \Xdot^{-1} \Xdot
\Ydot \Xdot  \Xdot^{-1} \prec \Ydot$$
yielding $i(X)\prec Y$. Similarly, by (2), we arrive at $Y\prec
i(X)$. Putting this together, we obtain the desired uniqueness
result. This shows that $\Ogamma$ is indeed an inverse semigroup.  Using
this, it is not hard to obtain the statement about the order. \hfill \qedsymbol

\medskip

Now, we combine this construction with the results of the first section
on groupoids associated to inverse semigroups.

\begin{definition}(a) The groupoid $\Gugamma\equiv G(\Ogamma)$ is called 
  the universal groupoid associated to $\Gamma$.\\
(b) The groupoid $\Gmgamma\equiv M(\Ogamma)$ is called the minimal
groupoid of $\Gamma$. 
\end{definition}

From the considerations of the first section, in particular, Proposition 
\ref{reduction}, we immediately infer the
following proposition. 

\begin{prop} $\Gmgamma= \Gugamma_{\Ogamma_{\rm min}^{(0)}}$. 
\end{prop}

Let us also note the following simple fact. 

\begin{prop} Let $\Gamma$ be an inverse semigroup with zero. Then
  $[B]\neq 0$ holds for every directed set $B$ with $0\notin B$. 
\end{prop}

\section{The groupoid $\Gugamma$}\label{ThegroupoidGu}
In this section we introduce a topology on $\Gugamma$ making it into a
topological $r$-discrete groupoid. This topology has a basis consisting
of compact $\Gugamma$-sets. We also show, that $\Gugamma$ with this
topology is actually isomorphic to the universal groupoid introduced by
Paterson in \cite{Pat1,Pat2}. 

\medskip

In the sequel we simply write $x$ instead of $[\{x\}]\in \Ogamma$ for
$x\in \Gamma$. In particular, we write $X\prec x$ instead of $X\prec
[\{x\}]$ for $X\in \Ogamma$. Note that we have $X= x X^{-1} X= X X^{-1} 
x $ for $X\prec x$. This will be used several times in the sequel. 
 For $x\in \Gamma$, we set $U_x\equiv \{X\in 
\Gugamma: X\prec x\}$.  For $x$, $x_1,\ldots, x_n\in \Gamma$ with
$x_1,\ldots, x_n\prec x$, we set
$$ U_{x; x_1, \ldots, x_n} \equiv U_x \cap U_{x_1}^c \cap\ldots \cap
U_{x_n}^c.$$
Here, $U_{x}^c$ is the complement of $U_x$ in $\Ogamma$. 
We will show that the family of these $ U_{x; x_1, \ldots, x_n} $  gives a
basis of a topology. To do so, we need the following proposition.

\begin{prop} \label{topology} For $X\in \Gugamma$ and $x_1, \ldots, x_n\prec x$ and
  $y_1,\ldots, y_m\prec y$ in $\Gamma$ with $X\in U_{x;x_1,\ldots, x_n}
  \cap U_{y;y_1,\ldots,y_m}$, there exist $z_1,\ldots, z_k\prec z$ with
  $z\prec x,y$ and 
  $ X\in U_{z;z_1,\ldots, z_k}\subset U_{x;x_1,\ldots, x_n}
  \cap U_{y;y_1,\ldots,y_m}$. 
\end{prop}
{\it Proof.} Let $p_j$ and $q_l$ in $\Gammanull $ be given with $x_j = x 
p_j$ and $y_l= y q_l$, $j=1, \ldots, n$, $l=1, \ldots, m$. By $X\in
U_x\cap U_y$, there exists $z\in \Gamma$ with $X\prec z\prec x,y$. Thus, 
there exist $p,q\in \Gammanull$ with $z=x p = y q = x p q = y p q$. Of
course, it suffices to show
$$X\in U_{z; z p_1,\ldots, z p_n, z q_1, \ldots, z
  q_m}\subset U_{x;x_1,\ldots, x_n} \cap U_{y;y_1,\ldots,y_m} .$$
It is straightforward to show that $X$ belongs to $U_{z; z p_1,\ldots, z p_n, z q_1, \ldots, z
  q_m}$. So, let us now show that $Y\in U_{z; z p_1,\ldots, z p_n, z q_1, \ldots, z
  q_m}$ belongs to $ U_{x;x_1,\ldots, x_n} \cap U_{y;y_1,\ldots,y_m}$ as 
well.  By $Y\prec z$ we have $Y\prec x$ and $Y\prec y$. Thus, it remains 
to show that $Y$ does neither  belong to $U_{x_j}$ nor to $U_{y_l}$ for
arbitrary  $j$ and $l$ as above. Assume $Y\prec x p_j$. By $Y \prec z$,
this gives the contradiction
$$ Y= Y Y^{-1} Y\prec x p_j z^{-1} x p_j= x p_j z^{-1} z p_j= x z^{-1} z 
p_j =z p_j,$$
where we used $z\prec x$ twice. Similarly, we show that $Y\prec y q_l$
cannot hold. The proposition follows. \hfill \qedsymbol

\medskip

The proposition implies that the family of all sets in $\Gugamma$ 
which are a union of sets of the form $U_{x;x_1,\ldots, x_n}$ is a
topology. 

\begin{definition} The topology $\Tau$ on $\Gugamma$ is the family of
  sets which are unions of sets of the form $U_{x;x_1,\ldots,
    x_n}$. 
\end{definition}

\begin{prop}\label{top} Inversion and multiplication in $\Gugamma$ are continuous
  with respect to $\Tau$. 
\end{prop}
{\it Proof.} The statement about inversion is obvious. To show that
multiplication is continuous, let $Z=X * Y\in U_{z;z_1,\ldots, z_n}$ be
  given. Let $p_j, q_j\in\Gammanull$ be given with $z_j=p_j z= z q_j$ for
  $j=1, \ldots, n$. There exist $x,y\in \Gamma$ with $X\prec x$, $Y\prec y$ and $x y\prec z$. As $X *  Y$ exists in $\Gugamma$, we
  have $X^{-1} X= Y Y^{-1}$ and we can assume w.l.o.g.  $x^{-1} x = y 
  y^{-1}$. Now, it is straightforward to show that $X\in U_{x;p_1 x,\ldots,
    p_n x}$ and $Y\in U_{y,y q_1,\ldots, y q_n}$. Thus, it remains to
  show that for $A\in U_{x;p_1 x,\ldots,
    p_n x}$  and $B\in U_{y,y q_1,\ldots, y q_n}$ the product $A * B$
  belongs to $U_{z;z_1,\ldots, z_n}$ (if it exists). Apparently, $A * B$ belongs to $U_{x
    y}\subset U_z$. Assume $A *  B \prec z_j$ for some $j$. Then, there
  exist $a,b\in \Gamma$ with $A\prec a$ and $B\prec b$ and $a b\prec
  z_j$. Again, as $AB$ exists in $\Gugamma$, we can assume
  w.l.o.g. $a^{-1} a= b b^{-1}$. Moreover, we can assume
  w.l.o.g. $a \prec x$ and $b\prec y$ as $A\in U_x$ and $B\in
  U_y$. This gives
\begin{eqnarray*}
a= a a^{-1} a &\prec&  a x^{-1} x= a a^{-1} a x^{-1} x= a b b^{-1}
x^{-1} x\\
&\prec & z_j b^{-1} x^{-1} x \prec z_j y^{-1} x^{-1} x  \prec p_j z
z^{-1} x\prec p_j x.
\end{eqnarray*}
This gives a contradiction, as $A$ does not belong to $U_{p_j x}$. The
proposition follows. \hfill \qedsymbol. 

\medskip

The proposition says that $\Gugamma$ with the topology $\Tau$ is a
topological groupoid. 
Let us now further investigate the topology. Even though the topology need not be Hausdorff, it has certain
separation properties. The following proposition shows in particular,
that the topology is $T_1$. Thus, a converging net cannot have more than one limit. 
\begin{prop}\label{Teins}(a) For arbitrary $X\neq Y\in \Gugamma$, there  exists $z\prec  x\in
  \Gamma$  with $X\in U_{x;z}$ and $Y\notin U_{x;z}$.\\
(b) The set $\Gugamma^{(0)}$ is closed in $\Gugamma$. 
\end{prop}
{\it Proof.} (a) Consider first  the case $Y\prec  X$ (and $X\neq Y$). Let
$x\in\Gamma$ with $X\prec x$ be given. Then, there exists an $y\prec x$
with  $Y\prec y$ and not $X\prec y$.  This gives $X\in U_{x;y}$ and
$Y\notin U_{x;y}$. On the other hand if $Y\prec X$ does not hold, then
there exists an  $x\in \Gamma$ with $X\prec x$ and not $Y\prec x$ and we 
infer $X\in U_x$ and $Y\notin U_x$. \\
(b) It suffices to show that, for  every converging net  $(P_i)$  in
$\Gugamma^{(0)}$, the limit $P$ belongs to $\Gugamma^{(0)}$
i.e. satisfies $P=P P^{-1}$. But this is immediate from (a) and continuity 
of multiplication. \hfill \qedsymbol

\begin{prop}\label{product} Let $p_1,\ldots, p_n\prec p\in \Gammanull$ be given. Let
  $x\in \Gamma$ with $p\prec x^{-1} x$ be given. Then $U_{xp;x
    p_1,\ldots, x p_n}= x U_{p;p_1,\ldots, p_n}$.
\end{prop}
{\it Proof.} This follows easily from $X= x X^{-1} X$ and the fact that
$X^{-1} X\prec q$ if and only if $x X^{-1} X\prec x q$ for $q\prec
x^{-1} x$ and $X\prec x$. \hfill  \qedsymbol

\medskip
 
Combining the  foregoing propositions, we infer the following
corollary. 

\begin{coro}\label{GSet} The maps $s:U_{x;x_1,\ldots,x_n}\longrightarrow U_{x^{-1} x; 
    x_1^{-1} x_1, \ldots, x_n^{-1} x_n}$, $X\mapsto X^{-1} X$ and
  $r:U_{x;x_1,\ldots,x_n}\longrightarrow U_{x x^{-1}; 
    x_1 x_1^{-1}, \ldots, x_n  x_n^{-1} }$, $X\mapsto X X^{-1} $ are
  homeomorphisms. 
\end{coro}
{\it Proof.} We only show the statement about $s$. The statement about
$r$ follows similarly. By the foregoing proposition, the map $s^{*}:U_{x^{-1} x; 
    x_1^{-1} x_1, \ldots, x_n^{-1} x_n} \longrightarrow
  U_{x;x_1,\ldots,x_n}$, $P\mapsto x P$ is surjective. By
$$ (x P)^{-1} x P= P x^{-1} x P= P P = P,$$
$s^{*}$ is injective as well. Moreover, we see that $s$ and $s^{*}$ are
inverse to each other and $s$ is therefore a  bijection. By Proposition
\ref{top}, the map $s$ is continuous.   Using
Proposition \ref{topology} and Proposition \ref{product}, one can also
infer that 
$s^{*}$ is continuous. \hfill \qedsymbol

\medskip

Recall that a groupoid is called $r$-discrete if its topology has a
basis of sets  on which $r$ and $s$ 
are homeomorphic. Thus, the foregoing corollary says that $\Gugamma$ is $r$-discrete.

\begin{prop}\label{compact} For  arbitrary $x_1,\ldots,x_n\prec x\in \Gamma$ the set
  $U_{x;x_1,\ldots,x_n}$ is compact. 
\end{prop}
{\it Proof.} By the foregoing Corollary, it suffices to consider
$U_{p;p_1,\ldots,p_n}$ with $p_1,\ldots, p_n\prec p\in \Gammanull$.  This, however, is just a reformulation of the well known properties of the maximal 
ideal space of the commutative Banach algebra $l^1(\Gammanull)$
(cf. \cite{Pat2} as well). We
include a short sketch for completeness.  Apparently, the map
$j:\Gugamma^{(0)}\longrightarrow \{0,1\}^{\Gammanull}$ with
$j(P)(q)=1$ if $P\prec q$ and $j(P)(q)=0$ otherwise, is
injective (cf. Lemma \ref{einfach} as well). Moreover, if $\{0,1\}$ carries the discrete topology and
$\{0,1\}^{\Gammanull}$ is  given the product topology, then 
the topology in $\Gugamma^{(0)}$ is easily seen to
be the topology induced by this injection. Thus, it remains to show that 
$j(\Gugamma^{(0)}$ is closed in $\{0,1\}^{\Gammanull}$. So, assume that
the net $(j(P_i))$ converges to $f\in \{0,1\}^{\Gammanull}$. Then, it is 
not hard to see that  $\{q\in \Gammanull: f(q) =1\}$ is a  directed
inverse order ideal and $f=j(P)$ with
$P=[\{ q\in \Gammanull   : f(q)=1\}]$. \hfill \qedsymbol

\medskip

Finally, we have the following proposition concerning the algebraic
properties of $x\mapsto U_x$. 

\begin{prop}\label{hom} The map $U: \Gamma \longrightarrow S(\Gugamma)$,
  $x\mapsto U_x$ is an injective homomorphism of inverse semigroups. 
\end{prop}
{\it Proof.} By Corollary \ref{GSet} the sets $U_x$ are indeed
$\Gugamma$-sets. Thus, $V$ maps into $S(\Gugamma)$. Apparently, $V$
preserves the involution. Thus, it only remains to show $U_{xy}= U_x
U_y$. The inclusion $\supset$ is obvious. Let now $Z\in U_{xy} $ be
given. By $Z\prec x y$, we have $x^{-1} Z\prec y$ and $Z y^{-1}
\prec x$ as well as  $Z y^{-1} x^{-1} = Z Z^{-1}$. This implies
$Z= Z Z^{-1} Z= Z y^{-1} x^{-1} Z= X Y$ with $X\equiv Z y^{-1} $ and
$Y\equiv x^{-1} Z$. It remains to show that $X$ and $Y$ are composable
in the sense of the groupoid $\Gugamma$ i.e. that $X^{-1} X= Y
Y^{-1}$. But this follows from 
$$X^{-1} X = y Z^{-1} Z y^{-1} = x^{-1} Z y^{-1} = x^{-1} Z Z^{-1} x = Y
Y^{-1},$$
where we used $Z y^{-1}
\prec x$ and $x^{-1} Z\prec y$. Injectivity is simple. \hfill  \qedsymbol. 

\medskip

We summarize our considerations  in the following theorem. 

\begin{theorem} The groupoid $\Gugamma$ is a topological groupoid with
  basis of topology given by the family of sets $U_{x;x_1,\ldots, x_n}$
  for  arbitrary $x_1,\ldots,x_n\prec x\in \Gamma$. These sets are compact
  $\Gugamma$-sets on which $r$  and $s$ are homeomorphisms.  The map $U:
  \Gamma\longrightarrow S(\Gugamma)$ is an injective homomorphism of
  inverse semigroups.  
\end{theorem}

Let us now consider the Hausdorff properties of $\Gugamma$. By the proof 
of Proposition \ref{compact}, its unit space is Hausdorff. However, in
general $\Gugamma$ will not be Hausdorff. We will show that, for $\Gamma$
satisfying condition (L), there  a simple alternative description of the
topology of $\Gugamma$. This will then  give that $\Gugamma$ is
Hausdorff if $\Gamma$ satisfies (L) (cf. Corollary \ref{Hausdorff}
below). 

Consider the map $j: \Gugamma\longrightarrow \{0,1\}^\Gamma$ with
$j(X)(x)=1$ if $X\prec x$ and $j(X) (x)=0$ otherwise. Let $\{0,1\}$
carry discrete topology and let 
$\{0,1\}^\Gamma$ be given the product topology.  We have the
following lemma.

\begin{lemma}\label{einfach} The map $j$ is injective. If $\Gamma$ satisfies (L), the
  topology induced on $\Gugamma$ from $ \{0,1\}^\Gamma$ agrees with $\Tau$.
\end{lemma}
{\it Proof.} It is not hard to show $X=[\{y: X\prec y\}]$, Thus, if
$X\neq Y$, then there exists w.l.o.g. an $x\in \Gamma$ with $X\prec x$
but not $Y\prec x$. This gives $j(X)(x)=1$ and $j(Y)(x)=0$ and
injectivity follows. 

To show that the induced topology agrees with $\Tau$, we have
to show that for arbitrary $X\in \Gugamma$ and $x_1,\ldots, x_n,
y_1,\ldots, y_m$ with  $X\in U_{x_1}\cap\ldots\cap U_{x_n} \cap
U_{y_1}^c\cap \ldots \cap U_{y_m}^c$, there exist $z_1,\ldots, z_k\prec
z$ with
\begin{equation} \label{stern}X\in U_{z;z_1,\ldots,z_k}\subset U_{x_1}\cap\ldots\cap U_{x_n} \cap
U_{y_1}^c\cap \ldots \cap U_{y_m}^c.
\end{equation}
By $X\in U_{x_1}\cap\ldots\cap U_{x_n}$, there exists an $x\in \Gamma$
with $X\prec x\prec x_1,\ldots, x_n$. By (L), we can then set $z\equiv
x_1\wedge \ldots \wedge x_n$. Apparently, we have $X\in U_z\subset
U_{x_1}\cap\ldots\cap U_{x_n}$. Similarly, we can define $z_j=z\wedge
y_j$ for every $j$ with $U_z\cap U_{y_j}\neq \emptyset$. Assume w.l.o.g.
that the set of these $j$ is given by $\{1,\ldots,k\}$. By construction
\eqref{stern} holds. \hfill  \qedsymbol. 

\begin{coro}\label{Hausdorff}
 If $\Gamma$ satisfies (L), then $\Gugamma$ is Hausdorff. 
\end{coro}

\begin{remark}{\rm In \cite{Pat2} it is shown that $\Gugamma$ is
    Hausdorff if $\Gamma$ is $E$-unitary. As $E$-unitary implies (L),
    the foregoing Corollary gives a strengthening of this
    result.} 
\end{remark}

We close this section with a  discussion  of the isomorphy between
$\Gugamma$ and the universal groupoid, $H_u (\Gamma)$ constructed by Paterson
\cite{Pat1,Pat2}. His construction  proceeds in three steps:

\begin{itemize}
\item A certain inverse  semigroup $\Gamma'$ containing $\Gamma$ is
  shown to act on the space $X$ of semicharacters s.t. $(X,\Gamma')$ is
  a localization. 
\item By Patersons extension of Kumjians theory of localizations
  \cite{Kum,Pat2}, there exist a groupoid
  $G(X,\Gamma')$ for this localization.
\item The groupoid $H_u(\Gamma)=G(X,\Gamma')$ can be expressed in terms
  of $X$ and $\Gamma$ only.
\end{itemize}

We refrain from discussing the theory of localizations here and just
give the description of $H_u (\Gamma)$ in terms of $X$ and
$\Gamma$ according to Theorem 4.3.1 of \cite{Pat2}. In our discussion we 
will identify the  space of semicharacters  used in \cite{Pat2} with 
$\Ogamma^{(0)}$  (cf. proof of Proposition \ref{compact} above and discussion in
Section 4.3 of \cite{Pat2}). Moreover, we will use the notation
introduced above. In particular,  the action of $\Gamma $ on $X$
will be written  accordingly. 
 Using these adoptions to our setting, the 
groupoid $H_u (\Gamma)$ can be described as follows:

It consists of equivalence classes  $[P,x]$, of pairs $(P,x)$ with  $P\in \Ogamma^{(0)}$, $x\in \Gamma$ with $P\prec x
x^{-1}$. Here, two pairs $(P,x)$ and $(\widetilde{P},\widetilde{x})$ are
identified  if $P=\widetilde{P}$ and there exists an $p\in \Gamma$ with
$P\prec p$ and $p x = p \widetilde{x}$.  
The involution is given by $[P,x]^*\equiv [x^{-1} P x, x^{-1}]$ 
and the multiplication is given by  $[P,x]* [x^{-1} P x, y]\equiv [P,x y 
]$. A basis of the topology is given by sets of the form $\{[P,x]: P\in
U_{p;p_1,\ldots,p_n}\}$. 

Given these reformulations of the Paterson construction, the proof of the following theorem is a simple exercise.

\begin{theorem} The map $J:
\Gugamma\longrightarrow H_u (\Gamma)$, $J(X)\equiv[X^{-1} X, x]$ with an 
arbitrary $x$ with $X\prec x$ is an isomorphism of topological groupoids 
with inverse map $K$ given  by $K:H_u (\Gamma) \longrightarrow\Gugamma$, 
$K([P,x])\equiv Px$. 
\end{theorem}

\section{The inverse semigroup $\Gammatilde$}\label{neu}
In this section we introduce and investigate a certain quotient of $\Gamma$, which we call $\Gammatilde$. The relevance of this quotient will become
apparent in the next sections when we deal with $\Gmgamma$. It 
will then turn out that $\Gammatilde$ and not $\Gamma$  is the appropriate
semigroup to phrase certain features of $\Gmgamma$. 

\medskip

\begin{definition} For $n\in \NN$ and $x,x_1,\ldots, x_n\in \Gamma$,  we
set $x<(x_1,\ldots,x_n)$ if for every $y\prec x$, $y\neq 0$, there
exists  $z\in \Gamma$, $z\neq 0$ and $j\in \{1,\ldots,n\}$  with
$z\prec y,x_j$.  If $n=1$, we write $x < x_1$ instead of   $x< (x_1)$.
\end{definition}

The relation $\gros$ is not an order. However, it can be shown to induce
an order on a certain quotient of $\Gamma$ by a standard procedure  in
the treatment of preorders. This is investigated next. 

The relation $x<>
y$ if and only if $x<y$ and $y<x$ can easily be seen to give an
equivalence relation on $\Gamma$. The quotient $\Gammatilde$ is then
defined by $\Gammatilde\equiv \Gamma/< >$. Let $\pi: \Gamma
\longrightarrow \Gammatilde$ be the canonical projection. 

Let us collect a few useful properties of $\gros$. 

\begin{prop} $x<y$ implies $x^{-1} < y^{-1}$ as well as   $x z< y z$ and $z x
< z y$. 
\end{prop} 
{\it Proof. } This is straightforward. \hfill \qedsymbol

\begin{prop} \label{eindeutig} There exists a unique  inverse semigroup structure on
$\Gammatilde$ making $\pi$ into a homomorphism of inverse semigroups.
The relation $x<y$ holds for $x,y\in \Gamma$ if and only $\pi(x)\prec
\pi(y)$. 
\end{prop}
{\it Proof.} The uniqueness statement is obvious. Let us now show
existence of the desired semigroup structure. Using the foregoing
proposition, we infer that the sets $\pi(\pi^{-1} (a)^{-1}  )$ resp. $\pi(\pi^{-1} (a) \pi^{-1} (b))$ contain exactly one element. Thus, we can define $ a^{-1}
$  resp. $b a$ by $\pi(\pi^{-1} (a)^{-1})$  resp. $\pi(\pi^{-1} (a) \pi^{-1} (b))$. Let
us now show that the inverse is unique. Let $x,y\in \Gamma$ and $a,b\in
\Gammatilde$ with $a=\pi(x)$ and $b=\pi(y)$ and $a b a= a$ and $b a b
=b$ be given. By $ a b a = a$, we have $x< x y x$ implying $x^{-1} = x^{-1} x
x^{-1} < x^{-1} x y x x^{-1} < y$. Similarly, we infer $y^{-1}  < x$ and $y^{-1} < > x$
follows. Thus, $\Gammatilde$ is indeed an inverse semigroup and $\pi$ is
an homomorphism of inverse semigroups. 

It remains to show the statement about the order. Let $x,y\in \Gamma$
with $x<y$ be given. We have to show $\pi(x) \pi(x^{-1} )  =  \pi(x)
\pi(y)^{-1}$ i.e. $x x^{-1} < > x y^{-1} $. By $x < y$ and the foregoing
proposition, we infer
\begin{equation}\label{tatutata} x x^{-1} < x y^{-1}, \;\:\;  x x^{-1} < y x^{-1}.
\end{equation}
Thus, it remains to show $x y^{-1} < x x^{-1}$. But this follows from 
$$ x y^{-1} = x x^{-1} x x^{-1} x y^{-1} < x x^{-1} y x^{-1} x y^{-1} \prec x x^{-1}.$$
Here, we used \eqref{tatutata} to obtain the estimate $<$.  Conversely,
assume $\pi(x)\prec \pi(y)$. This gives easily $x< y x^{-1} x \prec y$.
\hfill \qedsymbol
 
\medskip

Next, we study $\Gammatilde$ for $\Gamma$ satisfying (L). Our main tool
in this study is the following proposition. 

\begin{prop} \label{hilfe} Let $\Gamma$ satisfy (L). If $x,y,z\in
  \Gamma$ satisfy $0\neq z < x,y$, then $x\wedge y\wedge z$ exists and
  is not equal to zero.
\end{prop} 
{\it Proof.} By $z \prec z$ and $z < x$, we derive from (L) that $0\neq
z\wedge x$ exists. By $0\neq z\wedge x\prec z$ and $z < y$, we
infer, again by (L), that $0\neq z\wedge x\wedge y$ exists \hfill \qedsymbol

\medskip

We can now deduce two further properties of $\Gammatilde$. 

\begin{prop}(a)  If $\Gamma$ satisfies (L), then $\Gammatilde$ satisfies (L)
as well.\\
(b) Let $\Gamma$ be an inverse semigroup with zero satisfying
  (L). Then $\Gammatilde= (\Gammatilde)^{\widetilde{}}$
\end{prop}
{\it Proof.} (a) Let $0\neq c\prec a,b\in \Gammatilde$ be given. Choose
$x,y,z\in \Gamma$ with $\pi(x)=a$, $\pi(y)=b$ and $\pi(z)=c$. By
Proposition \ref{eindeutig}, we then have $0\neq z < x,y$. By the
foregoing proposition, $0\neq x\wedge y$ exists. By $x\wedge y
< x,y$ (even $x\wedge y\prec x,y$)  and Proposition \ref{eindeutig},  we
then have $\pi(x\wedge y)\prec \pi(x),\pi(y)$. 
Moreover, straightforward argument show that $z < x\wedge y$ holds
yielding $c=\pi(z)\prec \pi(x\wedge y)$.  Combining these estimates, we
infer $\pi(x\wedge y)= \pi(x)\wedge \pi(y)$. \\
(b)  It suffices to show $x< y$ whenever $\pi (x) < \pi (y)$ for
$x,y\in \Gamma$. So, assume $\pi(x) < \pi(y)$.  W.l.o.g. we can assume
$0\neq \pi(x)$. Let $0\neq z\prec x$ be given. Then,  we have $\pi(z)\prec \pi (x)$ and by $\pi (x) <  \pi (y)$,
there exists $r\in \Gamma$ with $0\neq \pi (r)\prec \pi
(z),\pi(y),\pi(x)$. This gives $0\neq r < z,y,x$ by Proposition
\ref{eindeutig}. By Proposition \ref{hilfe}, we then infer $0\neq
z\wedge y\wedge x\wedge r$  and $x < y$ follows. \hfill \qedsymbol

\medskip

For later use we also note the following proposition.

\begin{prop}\label{neueordnung} (a) The relation $\gros$ is transitive, i.e. for
  $x<(x_1,\ldots,x_k)$ with $x_j<(x_{j,1},\ldots, x_{j, n(j)})$,
  $j=1,\ldots, k$, $x_{j,l}\in \Gamma$ suitable,  the relation $x<(x_{1,1},\ldots,
  x_{1, n(1)},\ldots,x_{k,1},\ldots, x_{k,n(k)})$ holds.\\
(b) If $p<(p_1,\ldots,p_n)$ and $p \prec x^{-1} x$ for suitable
$p,p_1,\ldots, p_n\in\Gammanull$ and $x\in \Gamma$, then $x p x^{-1} <(x 
p_1 x^{-1},\ldots, x p_n x^{-1})$. 
\end{prop}
{\it Proof.} (a) This is straightforward.\\
(b) Let $0\neq q\prec x p x^{-1}$ be given. By $x p x^{-1}= x x^{-1} x p 
x^{-1} x x^{-1}$, this implies  $q = x x^{-1} q x x^{-1}$ yielding
$x^{-1} q x \neq 0$. Furthermore, we have $x^{-1} q x \prec x^{-1} x p
x^{-1} x =  p$. Thus, there exist $r\in
\Gammanull\setminus\{0\}$ and $j\in \{1,\ldots, n\}$ with $r 
\prec x^{-1} q  x $ and $r \prec p_j$. This implies $0\neq x
 r x^{-1}\prec x x^{-1} q x x^{-1}\prec q$ and $x
 r x^{-1} \prec x p_j x^{-1}$ finishing the proof of
 (b). \hfill \qedsymbol

\section{The groupoid $\Gmgamma$} \label{ThegroupoidGm}

By the general theory presented in Section \ref{Preliminaries} the
groupoid $\Gmgamma$ is a subgroupoid of $\Gugamma$ and in fact a
reduction in the set theoretical sense. Thus, it inherits the topology
from $\Gugamma$ and is a topological $r$-discrete groupoid. A basis of
the topology is given by the sets
$$ V_{x;x_1,\ldots, x_n} \equiv U_{x;x_1,\ldots,x_n} \cap \Gmgamma$$
for arbitrary $x_1,\ldots,x_n\prec x\in \Gamma$. 
However, it is not clear in general, whether the sets $V_{x;x_1,\ldots,
  x_n}$ are compact. Moreover, it is not clear whether $V_x$ is actually 
not empty.
So, we start by discussing conditions on non-emptyness and
compactness of $V_x$, $x\in\Gamma$. 

\begin{prop}\label{allesgut}
If $\Gamma$ contains a zero, then there exists for every $Y\in
\Ogamma$, $Y\neq 0$, an $X\in \Ogamma_{\min}$ with $X\prec Y$. In
particular, $V_x\neq \emptyset$ for every $x\neq 0$. 
\end{prop} 
{\it Proof.} Let $\Ydot$ be a representative of $Y$. By $Y\neq 0$, we have $0\notin \Ydot$. Consider the family
of directed sets containing $\Ydot$ but not containing $0$. The usual
inclusion gives a partial order on this family. Application of Zorns
Lemma, then gives a maximal element $B$ in this family. This element
does not contain zero and as $\Gamma$ contains a zero, we see $[B]\neq
0$. By construction $[B]$ is minimal and precedes   $ Y$. \hfill
\qedsymbol

\medskip

\begin{prop}\label{crucial}
(a) Let $\Gamma$ be  an inverse semigroup with zero satisfying $(L)$.
Then, the following are equivalent:

(i) $x<(x_1,\ldots,x_n)$.
(ii) $V_x\subset V_{x_1}\cup \ldots \cup V_{x_n}$. \\

In particular,  $V_x=V_y$ if and only if $x<y$ and $y<x$. \\
(b) For arbitrary  $\Gamma$ with zero (not necessarily satisfying (L)), the
equivalence of (i) and (ii) holds, whenever $x,x_1\ldots,x_n$ belong all
to $\Gammanull$. 
\end{prop}
{\it Proof.} (a) (i)$\Longrightarrow$(ii). Let $X\in V_x$ be given.
Then, $A\equiv \{y: y\prec x, X\prec y\}$ is a representative of $X$.
Set 
$A_j\equiv \{ 
y\wedge x_j : y\in A 
\;\mbox{ s.t. $0\neq y\wedge x_j$  exists }
           \}$. 
By $(i)$ and (L), there exists a $j$ with $A_j\prec
A$. This gives $[A_j]\prec X$. As $\Gamma$ has a zero, we have
$[A_j]\neq 0$ and by minimality of $X$, we infer $X=[A_j]$. As $[A_j]$
belongs to $V_{x_j}$, the statement (ii) is proven.\\
(ii) $\Longrightarrow$ (i). Let $y\prec x$, $y\neq 0$, be given. As $\Gamma$
has a zero, there exists by Proposition \ref{allesgut} an  $Y\in \Ogamma_{min}$ with $Y\prec y$.
This implies $Y\in V_y\subset V_x$.  By (ii), we infer $Y\in V_{x_j}$
i.e. $Y\prec x_j$  for a suitable $j$. Thus, $y$ and $x_j$ have a common
precessor not equal to zero. \\
(b) This follows easily from  the fact that existence of largest
precessors is always valid on $\Gammanull$. 
\hfill \qedsymbol

\medskip

We can  now study compactness properties of the $V_x$, $x\in \Gamma$.

\begin{prop}\label{equival}  The following are equivalent:\\
(i) For arbitrary $x_1,\ldots, x_n\prec x\in \Gamma$ the set
$V_{x;x_1,\ldots, x_n}$ is   compact.\\
(ii) The set $\Gmgamma^{(0)}$ is closed in $\Gugamma$.\\
(iii) $\Gmgamma$ is  a topological reduction of $\Gugamma$.
\end{prop}
{\it Proof.} The equivalence of (ii) and (iii) is immediate from
Proposition \ref{reduction}. The implication (ii) $\Longrightarrow$ (i)
is immediate from Proposition \ref{product} and Proposition \ref{compact}. Thus, it remains to show
(i) $\Longrightarrow$ (ii). Let $(P_i)$ be a net in $\Gmgamma^{(0)}$
converging in $ \Gugamma$ to $P\in \Gugamma$. By Proposition \ref{Teins} (b), $P$ belongs to $\Gugamma^{(0)}$. Thus, $P$ belongs to $U_p$ for a suitable $p\in
\Gammanull$. Then $P_i$ belongs to $V_p$ for large $i$. As $V_p$ is
compact,  $P=\lim P_i$ must then
belong to $V_p\subset \Gmgamma^{(0)}$ as well. \hfill \qedsymbol

\medskip

Of course, it might be difficult to decide whether one of the conditions 
of the proposition holds. Thus, let us now give a criterion, which
despite its simplicity can be checked for certain concrete semigroups
e.g those arising in the context of  tilings and graphs. 

A function $R: \Gamma\longrightarrow I$ with $I=[0,\infty)$ or
$I=[0,\infty]$  is called a
radiusfunction if it satisfies (R1) $R(x^{-1})=R(x)$, (R2) $R(x y) \geq
\min\{ R(x), R(y)\}$, (R3) $R(y)\leq R(x)$ for $x\prec y$. 
A radiusfunction $R$ on $\Gamma$,  gives rise to a
radiusfunction on $\Ogamma$, called $R$ again 
by $R(X)\equiv \sup \{R(x):  X\prec x \}$. A radius
function is called admissible if $R(X)=\infty$ if and only if $X\in
\Ogamma_{\rm min}$.

\begin{remark}{\rm  The definition of radiusfunction just says that
    $R$ is a dual prehomomorphism from $\Gamma$ into
    $(I,\wedge)$. Here, we set $x y \equiv x \wedge y \equiv
    \min\{x,y\}$ for $x,y\in I$.} 
\end{remark}

\begin{prop} \label{rad} If $R$ is an admissible and continuous radiusfunction on
  $\Gugamma$, then $V_{x;x_1,\ldots,x_n}$ is compact for arbitrary
  $x_1,\ldots, x_n\prec x\in \Gamma$. 
\end{prop}
{\it Proof.} By the foregoing proposition, it suffices to show that
$\Gmgamma^{(0)}$ is closed in $\Gugamma^{(0)}$. But this is immediate
from continuity of $R$ and $\Gmgamma^{(0)}=\Gugamma^{(0)}\cap \{X\in \Ogamma: R(X)=
\infty\}$. \hfill \qedsymbol

\begin{remark}{\rm (a) It is not hard to show that any radiusfunction must
    be lower semicontinuous.\\
(b) The radiusfunctions arising in the context of graphs or tilings are
admissible and  have 
the additional property that for $x\in \Gamma$ the maximum
$\pi_r(x)\equiv \max\{y: x\prec y, R(y)\geq r\}$ exists for arbitrary
$r\leq R(x)$. In this case it is possible to show that $\pi_r(X)\equiv
\pi_r(x)$ for arbitrary $x\in \Gamma, X\in \Ogamma$ with $X\prec x$ and
$R(x)\geq r$ is well defined and satisfies (1) $\pi_r(X)=\pi_r (Y)$ for
$X\prec Y$ and (2) $\pi_s(X)=\pi_s( \pi_r(X))$ for $s\leq r$ and
$R(X)\geq r$. Thus, in this case one can find a canonical representative 
$(\pi_n(X))_{n\in\NN}$ of $X\in \Ogamma_{\rm min}$. This can be used to
show that $\Gmgamma$ can be considered as  a kind of metric completion
of $\Ogamma$ w.r.t. $d(X,Y)\equiv \exp(- \sup\{r\geq 0:
  \pi_r(X)=\pi_r(Y)\})$.  In the context of tilings this has been investigated in \cite{Kel2} }
\end{remark}

Let us give another condition for closedness of $\Gmgamma^{(0)}$ in
$\Gugamma^{(0)}$. This condition is local in the sense that it can be
checked by only considering $\Gamma$ (and not $\Ogamma$). 

\begin{definition} The inverse semigroup $\Gamma$ is said to satisfy the
trapping condition (T), if $\Gamma$ contains a zero  and for every $p,q\in
\Gammanull$ with $q\prec p$  there exist $p_1,\dots, p_n\in \Gammanull$
with $p_j\prec p$, $j=1,\ldots, n$, and
\begin{itemize}
\item  $p < (p_1,\ldots, p_n, q)$.
\item  For every $j\in \{1, \ldots, n\}$ either $p_j\prec q$ or $p_j
q=0$. 
\end{itemize}
\end{definition}

\begin{prop}\label{trap} Let $\Gamma$ satisfy (T). Then $\Gmgamma^{(0)}$ is closed in
$\Gugamma^{(0)}$. 
\end{prop}
{\it Proof.} Let $(P_i)$ be a net in $\Gmgamma^{(0)}$ converging in $\Gugamma$ to $P$. As
$\Gugamma^{(0)}$ is closed in $\Gugamma$, the element $P$ belongs to
$\Gugamma^{(0)}$. Next, we show that $P$ is not zero. Assume the contrary.
As $\Gamma$ contains a zero, this implies $P\in U_0$ yielding the
contradiction  $0\neq P_i\in U_0$ for large $i$. 

So, it suffices to show
that every $Q\neq 0$, $Q\prec P$  agrees with $P$. Let such a
$Q$ be given and assume $P\neq Q$. Then, there exist $p,q\in \Gammanull$ with $q\prec p$ and
$Q\prec q$, $P\prec p$ but not $P\prec q$. Choose $p_1,\ldots, p_n$
according to (T) for $q\prec p$. Then, we have $V_p\subset V_{p_1}\cup
V_{p_n}\cup V_q$ by (T) and Proposition \ref{crucial} (b). Then,  it is not
hard to see that there exists a subnet $(P_k)$ of $(P_i)$  converging to
$P$ as well  and $(P_k)\subset V_{p_j}$ for a suitable $j$.  By $P_k\in
V_{p_j}\subset U_{p_j}$ and compactness of $U_{p_j}$, we infer $P\in
U_{p_j}$ i.e. $P\prec p_j$. There are two cases:\\
{\it Case 1.} $p_j\prec q$: In this case we arrive at the contradiction
$P\prec p_j\prec q$.\\
{\it Case 2.} $p_j q=0$ : In this case we have $P q= P p_j q = P 0=0$
contradicting $0\neq Q= Q q \prec P q$. \hfill \qedsymbol

\medskip

If $\Gamma$ satisfies (L) the topology of $\Gmgamma$ has a particularly
nice  basis.

\begin{lemma} \label{wichtig} If $\Gamma$ satisfies (L) and has a zero, then the family of sets $V_x$, $x\in \Gamma$, is a basis of the topology of $\Gmgamma$.
\end{lemma}
{\it Proof.} It suffices to show 
that for arbitrary   $X\in \Gmgamma$ and
$z_1,\ldots,z_n\prec z\in  \Gamma$ with $X\in V_{z;z_1,\ldots,z_n}$, we
have $X\in V_x\subset V_{z;z_1,\ldots,z_n}$ for a suitable $x$. Assume
the contrary. 
Thus,  there exists $X\in \Gugamma$ s.t.  for every $x$ with $X\prec x$ the set $V_x\cap (V_{z_1}\cup \ldots \cup
V_{z_n})$ is not empty. We therefore must have $X=[A_j]$
with  a suitable $j$ for $A_j\equiv \{x :X\prec x, \; V_x\cap V_{z_j} \neq
\emptyset\}$. Assume w.l.o.g. $j=1$. By (L), then the minimum $x\wedge
z_1$ exists for arbitrary $X\prec x$ and is not zero. Moreover, the
construction gives   $[\{x\wedge z_1: X\prec x\}]\prec X. $ and
$[\{x\wedge z_1: X\prec x\}]$  is not zero, as $\Gamma$ has a zero. By minimality of $X$, this gives $X=[\{x\wedge z_1: x\in
\Xdot\}]$ and the contradiction $X\prec z_1$ follows. \hfill \qedsymbol

\medskip

Let us now study the map  map $V: \Gamma \longrightarrow S(\Gmgamma)$,
$x\mapsto V_x$.

\begin{prop}\label{isovau} The map $V: \Gamma \longrightarrow S(\Gmgamma)$,
$x\mapsto V_x$ is an homomorphism of inverse semigroups. If $\Gamma$
satisfies (L) and contains a zero, $V(\Gamma)$ is canonically isomorphic 
to $\Gammatilde$ by $V_x\mapsto \pi(x)$.
\end{prop} 
{\it Proof. } The first statement can be shown with the same proof as 
 Proposition \ref{hom}. The second statement then follows from
Proposition \ref{crucial}. \hfill \qedsymbol

\medskip

\begin{remark}  
The proposition shows, in particular, that the map $V$ on $\Gamma$ is,
unlike $U$, not necessary injective. Nevertheless, it is still
possible to show that $\Gmgamma$ is isomorphic to$G_m( V(\Gamma))$,
whenever $\Gamma$ satisfies $(L)$.
\end{remark}

In this context, we also  have the following result. Recall that the ample
semigroup of a groupoid $G$ is the inverse  semigroup consisting of all compact open
$G$-sets.  A groupoid is called ample if this semigroup is a basis of
the topology. The result shows that the
construction of $\Gmgamma$ does not yield anything new if $\Gamma$ is
already (large part of) an ample semigroup of an ample groupoid.
\begin{theorem}\label{consistent} Let $\Gamma$ be a subsemigroup of the
  inverse semigroup of an ample Hausdorff groupoid $G$.  Assume that $\Gamma$ is closed 
  under intersections (which implies  (L)) and that $\Gamma$ is a basis
  of the topology of $G$. Then $\Gmgamma\simeq G$. 
\end{theorem}
{\it Proof. } This is the analogue in our setting to a result of  \cite{Kel2}. Thus, we only briefly sketch
the idea. To each point $g\in G$ we associate the set $A_g$
consisting of all $x\in \Gamma$ with $g\in x$. This set is directed
i.e. belongs to $\Ogamma$, as
$\Gamma$ is  closed under intersections. Using that $G$ is
Hausdorff, one easily sees that $[A_g]$ must be minimal i.e. belong to
$\Gmgamma$. Conversely, using finite intersection property of compact
sets, it is not hard to see that $g(X)\equiv \cap_{X\prec x} x$ is not
empty for every $X\in\Gmgamma$. By minimality of $X$ and again as $G$ is 
Hausdorff, the set $g(X)$ must then consist of only one point, which
 is denoted by $g(X)$. The maps $g\mapsto [A_g]$ and $X\mapsto
g(X)$ are groupoid homomorphisms and inverse to each other. \hfill \qedsymbol

\medskip

The considerations of this section suggest to distinguish the class of
inverse semigroups with zero which satisfy (L) and give rise to an ample 
groupoid $\Gmgamma$. Thus, we introduce the following definition. 

\begin{definition}\label{LC} The inverse semigroup $\Gamma$ is said to
  satisfy condition   (LC), if it contains a zero, satisfies (L) and $\Gmgamma^{(0)}$ is closed in
  $\Gugamma^{(0)}$. 
\end{definition}

\begin{remark}{\rm We will see that the inverse semigroup arising from
    the graphs considered in \cite{KPRR} and those arising from the
    tilings in \cite{Kel2} satisfy (LC).}
\end{remark}

Let us close this section with a  comparison to the corresponding
results of \cite{Kel2}. In \cite{Kel2} almost-groupoids are considered. An 
almost-groupoid is essentially an inverse semigroup with zero, whose
zero has been removed and whose multiplication has been restricted
accordingly. Thus, the inverse semigroups underlying the considerations
in \cite{Kel2} all have a zero.  The set of totally ordered sequences modulo the obvious
equivalence relation is  shown in \cite{Kel2} to be an almost-groupoid
whose set  of minimal elements is a 
groupoid and even a topological groupoid if equipped with the topology
generated by the $V_x$, $x\in \Gamma$ (in our notation). Call it $H_m (\Gamma)$. 

The considerations of this section extend the 
corresponding considerations of \cite{Kel} in some ways. 

First of all, the relationship between $H_m(\Gamma) $ and
$\Gugamma$ is made  explicit. More precisely, we show that $\Gmgamma$ is 
a subgroupoid and even a set theoretical reduction of $\Gugamma$. Now,
it can easily be seen that the
groupoid $H_m (\Gamma)$ agrees with $\Gmgamma$ as a set, but the
topology might be different. Here, Lemma \ref{wichtig} is important. It
shows that $\Gmgamma$ and $H_m (\Gamma)$ agree as topological groupoids
if $\Gamma$ satisfies (L).

Second of all we study the question   whether  $\Gmgamma$ has a basis
of compact open sets. In \cite{Kel2} this question is only addressed in
the tiling case. Here, we show that the existence of such a basis is essentially equivalent to 
$\Gmgamma$ being a topological reduction of $\Gugamma$. Moreover,   we
give two  simple sufficient  criteria on $\Gamma$ for this being the
case. Both criteria are met in both the tiling and graph case. 

Finally, as a minor point, let us remark that
out treatment is slightly more flexible as we use directed sets rather than totally ordered sequences.

\section{Open invariant subsets of $\Gmgamma^{(0)}$}\label{idealtheorie}
It is well known \cite{Ren} that  in arbitrary locally compact groupoids
$G$, each open invariant subset $U$ of $G^{(0)}$ gives rise to an ideal
in $\Cred (G)$ which is canonically  isomorphic to $\Cred (G_U)$. If $G$ is an essentially
principal groupoid (s. below for definition), then every ideal in
$\Cred(G)$ arises in this way. Thus, the investigation of open invariant 
subsets of $\Gmgamma^{(0)}$ is of primary importance. 

In this section we relate the open invariant subsets of $\Gmgamma^{(0)}$ to
certain order ideals in $\Gammanull$ (Lemma \ref{character}).  If
$\Gmgamma$ is essentially principal, this gives a complete
characterization of the ideals in $\Cred (\Gmgamma)$ (Theorem
\ref{ideal}). For $\Gamma$ satisfying (LC), we also give a necessary and sufficient condition on
$\Gamma$ for $\Gmgamma^{(0)}$ not to admit nontrivial invariant open sets
(Lemma \ref{minimal}).  This gives in particular a necessary and
sufficient condition for simplicity of $\Cred(\Gmgamma)$ whenever
$\Gmgamma$ is essentially principal. The semigroups we have in mind are
those satisfying (LC), even though some results of this section are
actually valid for more general inverse semigroups.

\medskip

We start with a discussion of invariance. Let $X\in \Gmgamma$ be
given,. Let $x\in \Gamma$ with $X\prec x$ be given. Then, we have $X= X
X^{-1} X= P x = x Q$ with  $P= X X^{-1},Q= X^{-1} X$ in $\Gmgamma^{(0)}$. This shows $X^{-1} X
= x^{-1} P x$ and  $X X^{-1}= x Q x^{-1}$.
These considerations easily imply the following proposition. 

\begin{prop} (a) For a subset $E$ of $\Gmgamma^{(0)}$ the following are
equivalent:\\
(i) $E $ is invariant.\\
(ii) For every $P\in E$ and $x\in \Gamma$ with $P \prec x x^{-1}  $, the
element $x^{-1} P x$ belongs to $E$.\\
(b) $\Gmgamma_P^P\equiv \{X: X X^{-1} = X^{-1} X=P\}= \{P\}$ if and only
if every $x\in \Gamma$ with  $x^{-1}
P x =P$  and $P\prec x x^{-1}$ satisfies $P \prec x$. 
\end{prop}

\begin{definition}
An element $P\in \Gmgamma^{(0)}$ is called aperiodic if $\Gmgamma_P^P=\{P\}$. 
\end{definition}

Using this definition and the above proposition, we can reformulate the
definition of (essentially) principality for $\Gmgamma$ given in
\cite{Ren} as follows:  $\Gmgamma$  is principal if and only if  every
$P\in \Gmgamma^{(0)}$ is aperiodic. $\Gmgamma$ is essentially principal if
and only if in every closed invariant set $F$ the set of aperiodic
points is dense.

\begin{definition}
(a) A subset $I$ of $\Gammanull$ is called $\gros$-closed if  $p\in
\Gammanull$ belongs to $I$ whenever   $p<(p_1,\ldots,p_n)$ for  $p_1,\ldots, p_n\in I$. \\
(b) A subset $I$ of $\Gammanull$ is called invariant if $x p  x^{-1}$
belongs to $I$ for every $p\in I$ and $x\in \Gamma$ with $p\prec x^{-1} x$
\end{definition}

Note that an $\gros$-closed set is in particular an order ideal as
$p\prec q$ implies $p<q$. 

\begin{prop} \label{achtvier}
(a) Let $I$ be an arbitrary subset of $\Gammanull$. Then $\Cl(I)\equiv
\{p: p<(p_1,\ldots,p_n)\;\:\mbox{for suitable}\:\; p_1,\ldots,p_n\in I\}$
is the smallest $\gros$-closed subset of $\Gammanull$ containing $I$. \\
(b) If $I$ is an invariant order ideal in $\Gammanull$ then $\Cl(I)$ is
 the smallest  $\gros$-closed invariant subset of $\Gammanull$ containing 
 $I$.
\end{prop}
{\it Proof.} (a)  It suffices to show that $\Cl(I)$ is
$\gros$-closed. This follows easily from Proposition \ref{neueordnung} (a).\\
(b) By (a), the set $\Cl(I)$ is $\gros$-closed. It remains to show
invariance. Let $p<(p_1,\ldots, p_n)$ with $p_1,\ldots, p_n\in I$ and $x\in
\Gamma$ with $p \prec  x^{-1} x $ be given. 
This gives, by Proposition \ref{neueordnung} (b),  $p=p p p < (p p_1
p,\ldots , p p_n p)= ( p_1 p,\ldots,  p_n p)$. Applying Proposition
\ref{neueordnung} once more, we infer

\begin{equation} \label{uuu} 
x p x^{-1} < (x p_1 p x^{-1} , \ldots, x  p_n p x^{-1} ).
\end{equation}
As $I$ is an order ideal, we have $p_j p\in   I$ for $j=1, \ldots, n$.
As $p_j p\prec p$ we have furthermore  $p_j p\prec x^{-1} x$ and by
invariance of $I$ this gives  $x p_j p x^{-1}\in I$ for $j=1, \ldots, n$.
Combining this with \eqref{uuu}, we conclude (b). \hfill\qedsymbol

\medskip

The proof of the following proposition is straightforward. 

\begin{prop}(a) The set of $\gros$-closed invariant subsets of
$\Gammanull$ with the usual inclusion  as partial order is a lattice with $I\vee
J\equiv \Cl(I\cup J)$ and $I\wedge J\equiv I\cap J$. \\
(b) The set of open invariant subsets of $\Gmgamma^{(0)}$ with the usual 
inclusion as 
order is a lattice with $U\vee V\equiv U\cup V$ and $U\wedge V\equiv
U\cap V$. 
\end{prop}

\begin{definition}(a) The lattice in part (a) of the foregoing
proposition will be denoted by $\Igamma$.\\
(b) The lattice in part (b) of the foregoing proposition will be denoted
by $\Ugamma$. 
\end{definition}

Next, we  prove the first key result of this section. 

\begin{lemma}\label{character} Let $\Gamma$ satisfy (LC).  For $V$ in
  $\Ugamma$ the set $S_i (V)\equiv \{q\in \Gammanull: V_q\subset  V\}$ belongs to $\Igamma$. For $I\in \Igamma$ the
set $S_u (I)\equiv \cup_{q\in I} V_q$ belongs to $\Ugamma$.  The maps $S_u:
\Igamma\longrightarrow \Ugamma$, $I\mapsto S_u (I)$, and $I:
\Ugamma\longrightarrow \Igamma$, $U\mapsto S_i (U)$, are lattice isomorphism
which are inverse to each other. 
\end{lemma}
{\it Proof.}  It is straightforward (and does not use any assumptions on
$\Gamma$) to show that   $S_u (I)$ belongs to $\Ugamma$. Moreover, using
(L), $0\in \Gamma$  and  Proposition \ref{crucial}, it is not hard to see
that $S_i (V)$ belongs to $\Igamma$. 

Let us now show that $S_i$ and $S_u$ are inverse to
each other, i. e. that 

$$ (1) \hspace{1ex} S_i (S_u (I)))=I\;\:\mbox{and}\:\; (2) \hspace{1ex} S_u (S_i (V ))=V.$$ \\
(1). By $S_i(S_u(I))=\{q: V_q\subset  \cup_{p\in I} V_p\}$, we have
$S_i(S_u (I))\supset I$. Let conversely, $q$ with $ V_q \subset \cup_{p\in I}
V_p$ be given. By compactness of $V_q$, we have $V_q\subset V_{p_1}\cup
\ldots \cup V_{p_n}$ for suitable $p_1,\ldots, p_n\in I$. By Proposition
\ref{crucial}, this gives $q<(p_1,\ldots, p_n)$. As $I$ is
$\gros$-closed, we infer $q\in I$ and the proof of (1) is finished. \\

(2). $S_u (S_i (V))= \cup_{q\in S_i (V)} V_q = \cup_{q: V_q\subset  V} V_q =
V$. Here, we used in the last equation that the $V_x$, $x\in \Gamma$,
give a basis of the topology of $\Gmgamma$ by (L).

\smallskip
 
Apparently, the maps $S_i$ and $S_u$ respect the order. So, it remains to
show that they respect $\vee$ and $\wedge$ as well. This will be shown
next. In fact, $S_i (U\wedge V)=S_i(U) \cap S_i(V)$ is immediate and
$S_u(I\wedge J)=S_u (I) \cap S_u (J)$ follows easily as $p\wedge q= p q$ 
exists for $p,q\in \Gammanull$. Thus, it remains to  show $S_u(I\vee J) =
S_u (I)\vee S_u (J)$ and $S_i ( U\vee V)= S_i (U) \vee S_i (V)$. 
We have 
$$S_u (I\vee J)\equiv S_u (\Cl(I\cup J))=\bigcup_{q\in \Cl(I\cup J)} V_q =
\bigcup_{q\in I} V_q \cup  \bigcup_{p\in J} V_p= S_u (I)\cup S_u (J),$$
where we used Proposition \ref{crucial} combined with Proposition \ref{achtvier} in the previous to the last
equality. Also, we have
\begin{eqnarray*}
S_i (U\vee V)&=& \{q: V_q\subset U\cup V\}\\
(V_q\;\:\mbox{compact})\;&=&\{q: V_q\subset V_{q_1}\cup \ldots
V_{q_n}\cup V_{p_1}\ldots V_{p_k}, V_{p_j}\subset U, V_{q_l}\subset
V\}\\
(\mbox{Prop. \ref{crucial}})\;\:&=& \Cl(\{q: V_q\subset U\} \cup \{ q:
V_q\subset V\} )\\
&=& S_i (U)\vee S_i (V). 
\end{eqnarray*}

This finishes the proof of the lemma. \hfill \qedsymbol

\medskip

The previous lemma characterizes the open invariant subsets of
$\Gmgamma^{(0)}$ in terms of invariant $\gros$-closed subsets of
$\Gammanull$. An important question is whether there actually exist
nontrivial invariant open subsets of $\Gmgamma$. This question is
answered in the following lemma.

\begin{lemma} \label{minimal} Let  $\Gamma$ satisfy (LC). Then the following are equivalent.\\
(i) There do not exist non trivial $\gros$-closed invariant subsets of
$\Gammanull$. \\
(ii) For every $p,q\in \Gammanull$, there exist $x_1,\ldots, x_n$ with 
$x_j^{-1} x_j \prec p$, $j=1, \ldots, n$,  and $q< (x_1 x_1^{-1}, \ldots,
x_n x_n^{-1})$.\\ 
(iii) $\Gmgamma$ is minimal, i.e. for every $P\in \Gmgamma^{(0)}$ the orbit
of $P$ is dense in $\Gmgamma$.

\end{lemma}
{\it Proof.} It is well known that an $r$-discrete topological
groupoid   $G$ is minimal if and only if there do not exist any non
trivial invariant open subsets of $G^{(0)}$. Thus, the equivalence of (iii)
and (i) is immediate from the foregoing lemma. 

It remains to show the equivalence of (i) and (ii). 
 Obviously, (i) is equivalent to the statement that any
non-empty $\gros$-closed invariant subset of $\Gammanull$ contains 
every  unit. This means that for every $p\in \Gammanull$, $p\neq 0$, 
the set 
$$I_p\equiv \Cl(\{ x r x^{-1}: r \prec p, r\prec  x^{-1}x    \} )$$
contains every  $q\in \Gammanull$. This is the case if and only if for every
$q\in \Gammanull$, there exist $y_1,\ldots, y_n$ and $r_1,\ldots, r_n\in
\Gammanull$ with $r_j\prec y_j^{-1} y_j,p$ and $q<( y_1 r_1
y_1^{-1},\ldots, y_n r_n y_n^{-1})$. But this is equivalent to (ii) with $x_j\equiv y_j r_j$, $j=1,
\ldots,n$, (resp. $r_j= x_j^{-1} x_j$, and $y_j= x_j$). \hfill \qedsymbol

\medskip

As mentioned in the introduction to this section we are interested in
reductions of $\Gmgamma$ to open invariant sets.  In our setting these reductions of
$\Gmgamma$ can directly be described in terms of certain subsemigroups
of $\Gamma$. This will be shown next.

It is not hard to see that, for each invariant $I\subset \Gammanull$ the 
set $\Gamma_I\equiv \{x: x x^{-1}\in I\}=\{x: x^{-1} x \in I\}$ with
multiplication and involution from $\Gamma$ is an inverse subsemigroup
of $\Gamma$. 

\begin{prop} \label{redu} Let $\Gamma$ be an inverse semigroup with zero satisfying
  (L). Let $I$ be an invariant order ideal in $\Gammanull$. Then
  $V(I)\equiv \cup_{q\in I} V_q$ is an invariant open subset of
  $\Gmgamma^{(0)}$. The canonical embedding $j: \Gamma_I\longrightarrow
  \Gamma$, $x\mapsto x$ induces an isomorphism $J: G_m
  (\Gamma_I)\longrightarrow \Gmgamma_{V(I)}$, $X\mapsto [ \{ j (y) : X \prec  y\}]$, 
  topological groupoids. 
\end{prop}
{\it Proof.} As in the proof of Lemma \ref{character} we infer that
$V(I)$ is open and invariant. Direct calculations show that  $J : G_m (\Gamma_I)\longrightarrow \Gmgamma$ and   $P:\Gmgamma_{V(I)}\longrightarrow G_m (\Gamma_I)$,
$X\mapsto [\{x\in \Gamma_I:  X\prec x\}]$ are  continuous groupoid
homomorphism which is inverse to each other. \hfill \qedsymbol

\medskip

This proposition allows one to identify $\Cred(G_m (\Gamma_I))$ with
  $\Cred( \Gmgamma_{V(I)})$ which in turn can canonically be considered as 
  an ideal in $\Cred(\Gmgamma)$ by the results of \cite{Ren} mentioned
  at the beginning of this section. Using this identification, we   can
  state   the results on the ideal structure of
$\Cred(\Gmgamma)$. We will denote  the lattice of ideals of $\Cred(\Gamma)$ by
${\mathcal{I}(\Cred(\Gamma))}$.

\begin{theorem}\label{ideal}
Let  $\Gamma$ satisfy (LC). Assume that
$\Gmgamma$ is essentially principal. Then the following holds.\\
(a) The map $J: \Igamma\longrightarrow {\mathcal{I}(\Cred(\Gamma))}$,
$J(I)\equiv \Cred (G_m (\Gamma_I))\subset \Cred(\Gmgamma)$ is a
bijection of lattices.\\
(b) The $C^{\ast}$-algebra $\Cred(\Gmgamma)$ is simple if and only if for
every $p,q\in \Gammanull$, there exist $x_1,\ldots, x_n$ with  $x_j^{-1}
x_j \prec p$, $j=1, \ldots, n$,  and $q< (x_1 x_1^{-1}, \ldots, x_n
x_n^{-1})$. 
\end{theorem}
{\it Proof.} (a) This follows from Lemma \ref{character} and the
corresponding results of chapter II, Section 4 of  \cite{Ren}
(cf. Corollary 4.9 of   \cite{Ren2} as well).\\
(b) This follows from (a) and Lemma \ref{minimal}.\hfill \qedsymbol

\section{Application to graphs} \label{Applicationgraph} 
In this section we present an inverse semigroup based approach to the
groupoids $G (\graph)$ associated to graphs $\graph$. Combined with
the results of the previous sections, this will provide semigroup
based proofs for some results of \cite{KPRR} concerning the structure
of the open invariant subsets of $G(\graph)^{(0)}$. Let us also
mention that the ideal theory of \cite{KPRR} has recently extended by
Paterson \cite{Pat3} to a non locally finite situation. He introduces
the universal groupoid associated to graph inverse semigroups and then
studies the graph groupoid as obtained by a reduction process.

\medskip

Let $\graph=(E,V,r,s)$ be a directed graph \cite{KPRR} with set of
edges $E$ and set of vertices $V$ and the range and source map
$r,s:E\longrightarrow V$. We assume that $r$ is onto and that
$s^{-1}(v)$ is not empty for each $v\in V$. Moreover, we assume that
the graph $\graph$ is row finite, i.e., that $s^{-1}(v)\subset E$ is
finite for all $v\in V$.

A  path $\alpha$  of length $|\alpha|=n\in \NN$  is a sequence
$\alpha=(\alpha_1,\ldots,\alpha_n)$ of edges $\alpha_1,\ldots,
\alpha_n$ in $E$ with $s(\alpha_{j+1})=r(\alpha_j)$,
$j=1,\ldots,n-1$. For such an $\alpha$ we set $s(\alpha)\equiv
s(\alpha_1)$ and $r(\alpha)\equiv \alpha_k$.  A path of length $0$ is just a
vertex and will also be called a degenerate path.  For such a path $v$, we set $r(v)\equiv v$ and $s(v)\equiv
v$. 

The set of all paths of finite length is denoted by $F(\graph)$. The set
of all infinite paths $\alpha=(\alpha_1,\ldots)$ is denoted by
$P(\graph)$.  The concatenation $\alpha \mu$ of two finite paths $\alpha$ and 
$\mu$ with $r(\alpha)=s(\mu)$  is defined in the obvious way, i.e. by 

$$\alpha \mu\equiv 
\left\{ 
\begin{array}{r@{\quad:\quad}l}
(\alpha_1,\ldots, \alpha_{|\alpha|},\mu_1,\ldots, \mu_{|\mu|}) &
|\alpha|,|\mu|>0\\
\alpha & |\mu|=0\\
\mu & |\alpha|=0
\end{array}
\right.
$$
By a slight abuse of language we write $\alpha \prec \beta$ if
$\beta=\alpha \mu$. 

We will now introduce an inverse semigroup associated to $\graph$. This
inverse semigroup is slightly more general than the inverse semigroup
discussed  e.g. in \cite{Pat2} in that we also allow for paths of
lengths zero. \\
Let the set  $\Gamma\equiv \Gamma(\graph)$ be given by $\Gamma\equiv \{(\alpha,\beta)\in
F(\graph)\times F(\graph)\,:\, r(\alpha)=r(\beta)\} \cup \{0\}$. Let us
now define a multiplication on $\Gamma$:  For a
pair $((\alpha,\beta),(\gamma,\delta))$ with $\gamma=\beta \mu$ with a
suitable (possibly degenerate) $\mu$, we define the product by 
$$(\alpha,\beta)(\gamma,\delta)=(\alpha,\beta)(\beta \mu,\delta)\equiv
(\alpha \mu, \delta).$$

For a pair $((\alpha,\beta),(\gamma,\delta))$ with $\beta= \gamma \mu$ with
a suitable (possibly degenerate)  $\mu$ we define the product by 
$$(\alpha,\beta)(\gamma,\delta)=(\alpha,\gamma \mu)(\gamma,\delta)\equiv (\alpha,\delta \mu).$$

In all other cases we define the product to be $0$. 
It is then not hard to show that $\Gamma$ is indeed an inverse semigroup, where the inverse of  $(\alpha,\beta)$ is given by $(\alpha,\beta)^{-1}\equiv (\beta,\alpha)$.  
Thus, $\Gamma$ gives rise to a groupoid $\Gmgamma$.
Now,  the function $R:\Gamma\longrightarrow [0,\infty)$ given by
$$R(\alpha,\beta) \equiv 
\left\{ 
\begin{array}{r@{\quad:\quad}l}
0 & \alpha_{|\alpha|}\neq\beta_{|\beta|}\\
\sup\{j\in \NN_0: \alpha_{|\alpha|-i}= \beta_{|\beta|-i}, i=0, \ldots, j \}
  & otherwise
\end{array}
\right.
$$
can easily be seen to be a radiusfunction in the sense of Section
\ref{ThegroupoidGm}. Moreover, it is possible to show that $R$ is
admissible and continuous. Thus, $\Gmgamma$ is a groupoid with a basis consisting of
compact sets. We refrain from giving details, but rather show the
connection between $\Gmgamma$ and the graph groupoids which were
introduced in \cite{KPRR}. 

To do so, we start with a simple proposition giving an understanding of
the relation $\prec$ in the case at hand.

\begin{prop} The relation $(\gamma,\delta)\prec  (\alpha,\beta)$ holds if and only if there exists a (possibly degenerate) $\mu\in F(\graph)$ with $\gamma=\alpha \mu$ and $\delta=\beta \mu$. 
\end{prop}
{\it Proof.} This is straightforward.\hfill
\qedsymbol. 

\medskip

From this proposition we immediately infer the following interesting fact 
concerning the order structure of $\Gamma$.

\begin{prop} Let $x,y\in \Gamma$ with a common precessor  be given. Then
  either $x\prec y$ or $y\prec x$. 
\end{prop}

Thus, to every $X\in \Ogamma$, we can  find $(\alpha,\beta)\in
\Gamma$, $I=[0,a]\subset \ZZ$, $a\in  \NN\cup\{\infty\}$,  and edges $e_n$, $n\in I$, s.t. $\{(\alpha e_1\ldots e_n, \beta e_1\ldots e_n): n\in I\}$ is a
representative of $X$. 

Putting this together, we see that minimal elements in $\Ogamma$ can be
identified with double paths of infinite length which agree from a certain point on. But
this is exactly the way the graph groupoid in \cite{KPRR} is
constructed. Let us sketch the construction and give a precise proof of
the isomorphy. 
Two paths $x,y\in P(\graph)$ are called equivalent with lag $k\in \ZZ$, written as $x\sim_k y$, if there exists an $N(x,y)\in \NN$ s.t. $x_i=y_{i+k}$ for all $i\geq N(x,y)$. Let
$$G(\graph)\equiv\{(x,k,z)\in P(\graph)\times \ZZ\times P(\graph)\,:\,x\sim_k y\}.$$
Let the set of composable pairs of $G(\graph)$ consist of all pairs $((x,k,y_1),(y_2,l,z))$ with $y_1=y_2$. For such a pair define the multiplication by
$$ (x,k,y_1)(y_2,l,z)\equiv (x,k+l,z).$$
An inverse map on $G(\graph)$ is given by letting
$$ (x,k,y)^{-1}\equiv (y,-k,x).$$
Then $G(\graph)$ together with this multiplication and inverse map is a groupoid. To make $G(\graph)$ into a topological groupoid one introduces the sets
$$Z(\alpha,\beta)\equiv\{(x,k,y)\in G(\graph)\,:\, \alpha\prec x,
\beta\prec y, k=|\beta|-|\alpha|, x_i=y_{i+k}, i\geq |\alpha|\},$$
where $(\alpha,\beta)$ is an arbitrary element in $F(\graph)\times
F(\graph)$ with $r(\alpha)=r(\beta)$. These sets form a basis for a
locally compact topology on $G(\graph)$. Each set $Z(\alpha,\beta)$ is a
compact and open $G(\graph)$-set. It is not hard to show that the system
$\mathcal{Z}(\graph)$ of all these sets  is in fact an inverse semigroup.

\begin{theorem}\label{graphiso} 
The map $ j: G(\graph)\longrightarrow G_m (\Gamma(\graph))$,
$j((x,k,y)) :=[\{ (x_1,\ldots, x_n, y_1,\ldots, y_{n+k}) : n\geq
N(x,y)\}]$, is an isomorphism of topological groupoids.
\end{theorem}
{\it Proof.}  By definition of $G (\graph)$ and the previous
proposition, the set $\{ (x_1,\ldots, x_n, y_1,\ldots, y_{n+k}) :
n\geq N(x,y)\}$ is directed and contains elements of arbitrary
lengths. Thus, $j(X)$ is minimal.  Direct calculations show that $j$
is a groupoid homomorphism. Now, for each $X\in G_m (\Gamma(\graph)) $
we can find an $(\alpha,\beta)\in \Gamma(\graph)$ with $X \prec
(\alpha,\beta)$. As $X$ is minimal and $s^{-1} (v) \neq \emptyset$ for
each vertex $v$, there exists for each $n\in \NN$ an path $\mu_n$ of
length $n$ with $ X\prec (\alpha \mu_n, \beta \mu_n)$. By the previous
proposition
\begin{equation}\label{ast} 
(\alpha \mu_{n+1}, \beta \mu_{n+1}) \prec (\alpha \mu_n, \beta \mu_n)
\end{equation}
for every $n\in \NN$.  Then, $\{ (\alpha \mu_n, \beta \mu_n) : n\in
\NN\}$ is directed and $X':= [\{ (\alpha \mu_n, \beta \mu_n) : n\in
\NN\} ]$ is minimal, as it contains paths of arbitrary lengths. As by
construction $X\prec X'$, we infer $X=X'$. Moreover, by \eqref{ast},
we can form the ``limit'' $x$ of the paths $\alpha \mu_n$, the ``limit''
$y$ of the paths $\beta \mu_n$ and define the map $h : G_m
(\Gamma(\graph))\longrightarrow G(\graph)$ via
$$ h(X) :=(x,k,y)$$ with $k:=|\alpha|- |\beta|$. By construction, $h$
is inverse to $j$. Moreover, it is not too hard to see that $h$ is a
homomorphism of groupoids. Thus, $j$ is an isomorphism of groupoids.

Direct calculation show that $j(Z(\alpha,\beta)= V_{(\alpha,\beta)}$
and we see that $j$ and $h$ are continuous.  \hfill $\Box$

\medskip

This theorem allows us to apply the theory of the preceeding sections
to the study of graph groupoids. In particular, we can rephrase the
ideal theory of \cite{KPRR} (viz the characterization of open
invariant subsets of $G(\graph)^{(0)}$ ) in terms of inverse
semigroups using Section \ref{idealtheorie}. This is done next.

Recall the following definitions from \cite{KPRR}.  For vertices $v,w\in 
V$ we write $v \geq w$ if  there exists a path in $P$ from $v$ to $w$. 
A subset $H$ of $V$
is called hereditary if $v\in H$ and $v\geq w$ implies $w\in H$ and it
is called saturated if
$$\mbox{ [$r(e)\in H$ for all $e\in E$ with $s(e)=v$ ] implies $v\in
  H$. }$$ 
The set  of hereditary and saturated subsets of $V$ is a
lattice under the operation of intersection of sets and union followed by 
saturation. 
Now, we have the following lemma. 
\begin{lemma}\label{it} The map $I\mapsto \{r(p): p\in I\}$ is an isomorphism
  between the lattice of invariant $\gros$-closed ideals in $\Gammanull$ 
  and the lattice of hereditary saturated subsets of $V$.  The inverse
  is given by $H\mapsto \{p\in  \Gammanull: r(p)\in H\}$. 
\end{lemma}
{\it Proof.} This follows by direct arguments. \hfill  \qedsymbol

\section{Application to tilings}\label{tilings}
As mentioned in the introduction, our study is motivated by work of
Kellendonk \cite{Kel,Kel2} introducing inverse semigroups in the
context of tilings, see \cite{KL,KL2} for recent work on this.

Here, we shortly discuss how the groupoid arises from the inverse
semigroup in this context. This follows essentially \cite{Kel} (with
the slight variation that we work with directed sets rather than
directed sequences). We then, apply the general theory developed above
to describe the ideal structure of $\Cred ( \Gmgamma)$ for $\Gamma$
arising from aperiodic tilings. While this is essentially, known it
serves as a good example for our theory. Moreover, it underlines the
structural similarities between tilings and graphs.

\medskip

A tiling in $\RR^d$ is a (countable) cover $T$ of $\RR^d$ by compact
sets which are homeomorphic to the unit ball in $\RR^d$ and which
overlap at most at their boundaries \cite{GS}. The elements of $T$ are
called tiles. A pattern $P$ in $T$ is a finite subset of $T$. For
patterns $P$  and tilings $T$ and $x\in \RR^d$, we define $P+x $ and $T+ 
x$ in the obvious way. The set of
all  patterns which belong to     $T+x $ for some  $x\in \RR^d$,  will be denoted by $P(T)$.
All patterns will be assumed to be patterns in $P(T)$ if not stated
otherwise. 

A doubly pointed pattern $(a,P,b)$ (over $T$)  consists of a pattern
$P\in P(T)$ together with two tiles $a,b\in P$. We say that $(a,P,b)$ is
contained in $(c,Q,d)$, written as $(a,P,b)\subset(c,Q,d)$, if $a=c,$
$b=d$ and $P\subset Q$. On the set of doubly pointed patterns over $T$
we introduce an equivalence relation by defining $(a,P,b)\sim (c,Q,d)$
if and only if there exists an $r\in \RR^d$ s.t. $c=a+r$, $d=b+r$ and
$Q=P+r$.  The class of $(a,P,b)$ will be denoted by $\overline{(a,P,b)
  }$.Obviously, the relation $\subset$ can be extended to these
classes. 

Similarly, one can introduce an equivalence relation on the set of all
patterns in $P(T)$.  Denote the class of the pattern $P$ up to
translation by $\overline{P}$ and denote the set of all classes of
patterns in $P(T)$  by $\overline{P(T)}$.   
Following \cite{Kel,Pat2}, we will assume  two finite type  conditions:
\begin{itemize}
\item[(i)] $d_{max}\equiv \sup\{\diam(A)\,:\, A\in T\}<\infty$.
\item[(ii)] The  set $\{\overline{P}\in \overline{P(T)}\,:\, \diam(\cup_{t\in P } t)\leq R\}$ is finite for every $R$.
\end{itemize}
Here, $\diam(A)$ denotes the diameter of $A$. Note that these conditions
imply in particular that there only finitely many different tiles up to
translation. As each tile is homeomorphic to the unit ball, this implies
in particular that there is a minimal volume $V_{min}>0$ among the
volumes of the tiles.

Following \cite{Kel}, we make $\Gamma\equiv\{\overline{(a,P,b)
  }\,;\,P\in P(T), a,b\in P\}\cup \{0\}$ into an inverse semigroup in the following way. A pair $(E,F)\in \Gamma\times \Gamma$ is said to be composable if there exists a doubly pointed pattern class $ G$  and  representatives $(a,P,b) $ of $E$, $(c,Q,d)$ of $F$  and $(r,R,s)$ of $G$  together with a tile $t\subset R$  with\\
$$(a,P,b)\subset (r,R,t)\hspace{4ex}\mbox{and}\hspace{4ex} (c,Q,d)\subset(t,R,s).$$
Let $H$ be the smallest w.r.t. $\subset$ doubly pointed pattern class
with this property.   It is not hard to see that  $E F\equiv H$ is well
defined. Now, we can define a multiplication on $\Gamma$ by $E F=H$ if
$E,F$ are  composable and by $E F= 0$ otherwise.

It can be shown that $\Gamma$ with the above multiplication is
indeed an inverse semigroup with inverse map given by
$\overline{(a,P,b)}^{-1}\equiv  \overline{(b,P,a)}$.  

Moreover, the relation ``$\prec$''  induced from the almost-groupoid agrees with the
relation ``$\subset$'' defined above \cite{Kel}, i.e. the following is
valid:
\begin{prop} For $x,y\in  \Gamma$ the relation $x\prec y$ holds if and
  only if there exist representatives $(a,P,b)$ of $x$ and $(a,Q,b)$ of
  $y$ with $P\supset Q$.
\end{prop}

Of course, $\Gamma$ gives now rise to a groupoid $\Gmgamma$. This
groupoid can easily be identified with the groupoid $G(T)$ defined as
follows \cite{Kel}: 
Let $\mathcal{T}=\mathcal{T}(T)$ be the set of all tilings $S$ of $\RR^d$ with      $\overline{P(S)}\subset \overline{P(T)}$. Let $G(\mathcal{T})$ denote the set of
 all     equivalence classes of doubly pointed tilings of
 $\mathcal{T}$. Here, a 
 doubly pointed tiling and the equivalence relation are defined by just
 replacing the pattern $P$ in the corresponding definitions above by a
 tiling $S\in \mathcal{T}$.  The set $G(\mathcal{T})$ has a groupoid
 structure. Two elements $E,F$ are composable if there exists
 representatives 
$(a,S,b)$ of $E$ and $(c,R,d)$ of $F$  with $S=R$ and $b=c$.  In this case one defines $E F\equiv \overline{(a,S,d)}$. This is well defined. The topology on $G(T)$ is generated by the sets
$$V(\overline{(a,P,b)})=\{E\in G(T)\,:\, \overline{(a,P,b)}\subset E\}.$$
These sets are in fact compact,
 open $G(\mathcal{T})$-sets  forming a basis of the
topology.  As in \cite{Kel2}, one can then show that $\Gmgamma=G(\mathcal{T})$. 

\medskip

$\Gamma$ admits a complete radius function $R$,
where $R(\overline{(a,P,b)})$ is defined by $R(\overline{(a,P,b)}\equiv
[dist(\partial P, \{a,b\})]$. Here, $\partial P$ is the boundary of
$P$. 

\medskip

Let us now study the structure of open invariant sets in $\Gmgamma$. A
subset $S$ of $\overline{P(T)}$ is called saturated if $\overline{P}\in
\overline{P(T)}$ and $\overline{Q}\in S$ with $\overline{Q}\subset \overline{P}$  implies
$\overline{P}\in S$. A subset  $S$ of $\overline{P(T)}$ is called
hereditary if  $\overline{P}$ belongs to $S$, whenever there exist $
\overline{P_1},\ldots, \overline{P_n}$ in $S$ satisfying the following
condition:
\begin{itemize}
\item For every pattern  $Q$ with $R(Q)$  large enough and
  $Q\supset P$, there exists $j\in \{1,\ldots,n\}$ with $Q\supset P_j$.
\end{itemize}

Then, we can easily infer the following lemma.

\begin{lemma} The map $I \mapsto \{ \overline{P}: \overline{(a,P,b)}\in
  I\}$ is an   isomorphism of the lattice of invariant $\gros$-closed
  subsets of $\Gammanull$ and the lattice of saturated hereditary subsets of
  $\overline{P(T)}$. 
\end{lemma}

It remains to study principality of $\Gmgamma$. Here, we have a very
simple and well-known condition. Recall, that a tiling $S$ is called
periodic if there exists an $x\in \RR^d$ with  $S+ x = S$. Now,
$\mathcal{T}$ is called aperiodic if it does not contain a periodic
tiling. 

\begin{lemma} $\Gmgamma$ is principal if and only if $\mathcal{T}$ is
  aperiodic. 
\end{lemma}
{\it Proof.} $\Gmgamma$ is principal, if every $P$ in $\Gmgamma^{(0)}$
is aperiodic in the sense of Section \ref{idealtheorie}. But this can
easily be seen to be equivalent to $\mathcal{T}$ being aperiodic in the
sense given above. \hfill \qedsymbol

\medskip

{\it Acknowledgements.}  Very special thanks are due to J. Kellendonk
and M. Lawson for stimulating discussions and comments on inverse
semigroups and tilings.  The author would also like to thank
A. Paterson, who commented on an earlier version of this paper.

\end{document}